\RequirePackage{fix-cm}
\documentclass[smallextended,numbook,runningheads]{svjour3}     
\smartqed  

\usepackage{graphicx}
\usepackage{amsmath}
\usepackage{epstopdf} 
\usepackage{bbding}
\usepackage{mathptmx}
\usepackage{algorithm}
\usepackage{algorithmic}
\usepackage{color}

\begin{document}

\title{Local discontinuous Galerkin methods for fractional ordinary differential equations
\thanks{Supported by NSFC 11271173, NSF DMS-1115416, and OSD/AFOSR FA9550-09-1-0613.
}}


\author{Weihua Deng        \and
        Jan S. Hesthaven
}


\institute{Weihua Deng.  E-mail: dengwh@lzu.edu.cn \at
              School of Mathematics and Statistics, Lanzhou University, Lanzhou 730000, P.R. China \\
           \and    Jan S. Hesthaven.  E-mail: Jan.Hesthaven@epfl.ch         \at
           {EPFL-SB-MATHICSE-MCSS, \'{E}cole Polytechnique F\'{e}d\'{e}rale de Lausanne, CH-1015 Lausanne,
Switzerland}}

\date{Received: date / Accepted: date}

\maketitle

\begin{abstract}
This paper discusses the upwinded  local discontinuous Galerkin
methods for the one-term/multi-term fractional ordinary differential
equations (FODEs).
 The natural upwind choice of the numerical fluxes for the initial value problem for FODEs ensures  stability of the methods.
 The solution can be computed element by element with optimal order of convergence $k+1$ in the $L^2$ norm and  superconvergence of
 order $k+1+\min\{k,\alpha\}$ at the downwind point of each element. Here $k$ is the degree of the approximation polynomial used in
 an element and $\alpha$ ($\alpha\in (0,1]$) represents the order of the one-term FODEs. A generalization of this includes problems with
 classic $m$'th-term FODEs, yielding superconvergence order at downwind point as $k+1+\min\{k,\max\{\alpha,m\}\}$. The underlying mechanism of
 the superconvergence is discussed and the analysis confirmed through examples, including a discussion of how to use the scheme as an efficient
  way to evaluate  the generalized Mittag-Leffler function and solutions to more generalized FODE's.
\keywords{fractional ordinary differential equation \and local
discontinuous Galerkin methods \and downwind points \and
superconvergence \and generalized Mittag-Leffler function}
\subclass{34A08  \and 65L60 \and 33E12}
\end{abstract}

\section{Introduction}

Fractional calculus is the generalization of classical calculus and
with a history parallel to that of classical calculus. It was
studied by many of the same mathematicians who contributed
significantly to the development of classical calculus. In recent
years, it has seen a growing use as an appropriate model for non
classic phenomena in physical, engineering, and biology fields
\cite{Metzler:00}.

To gain an understanding of the basic elements of fractional
calculus, we recall
\begin{equation}\label{DefFracInt1}
\int_a^t d \tau_n \int_a^{\tau_n} d \tau_{n-1} \cdots
\int_a^{\tau_2}x(\tau_1) d \tau_1=\frac{1}{(n-1)!}\int_a^t
(t-\tau)^{n-1}x(\tau)d\tau,\,\,t>a.
\end{equation}
A slight rewriting yields
\begin{equation}\label{DefFracInt2}
{_aD_t^{-n}}x(t)=\frac{1}{\Gamma(n)}\int_a^t
(t-\tau)^{n-1}x(\tau)d\tau,\,\,t>a.
\end{equation}
Clearly, Eq. (\ref{DefFracInt2}) still formally makes sense if  $n$
is replaced by $\alpha$ ($\in {R}^{+}$) leading to the definition of
the fractional integral as
\begin{equation}\label{DefFracInt3}
{_aD_t^{-\alpha}}x(t)=\frac{1}{\Gamma(\alpha)}\int_a^t
(t-\tau)^{\alpha-1}x(\tau)d\tau,\,\,t>a,\,\,\alpha \in {R}^+.
\end{equation}
Here $a \in {R}$. This definition of the fractional integral is
natural and simple, and inherits many good mathematical properties,
e.g., the semigroup property. Moreover $-\alpha$ `connects'
$-\infty$ to $0$ in a natural way.

With the derivative as the inverse operation of the integral, the
definition of the fractional derivative naturally arises by
combining the definitions of the classical derivative and the
fractional integral.  In this way, the fractional derivative
operator $D^\alpha$ can be defined as $D^{\lceil \alpha
\rceil}D^{-(\lceil \alpha \rceil-\alpha)}$ or $D^{-(\lceil \alpha
\rceil-\alpha)}D^{\lceil \alpha \rceil}$. Mathematically, the first
one can be regarded as preferred since $D^nD^{-n}=I$ while
$D^{-n}D^n=I+\cdots$, which involves additional
 information about the function  at the left end point. When considering a fractional derivative in the temporal direction, the second formulation is more popular due to its convenience when specifying the initial condition in a classic sense.
 These two definitions of derivatives are referred to as the  Riemann-Liouville derivative and Caputo derivative \cite{Podlubny:99}, respectively, with the
Riemann-Liouville derivative being defined as
\cite{Butzer:00,Podlubny:99}
\begin{equation}\label{RLDeri}
{_aD_t^\alpha}x(t)=\frac{1}{\Gamma(n-\alpha)} \frac{d^n}{dt^n}
\int_a^t (t-\tau)^{n-\alpha-1}x(\tau)d \tau,~~ t>a,\,\, \alpha \in
[n-1,n);
\end{equation}
and the Caputo derivative recovered as
\begin{equation}\label{CDeri}
{_a^CD_t^\alpha}x(t)=\frac{1}{\Gamma(n-\alpha)} \int_a^t
(t-\tau)^{n-\alpha-1} \frac{d^nx(\tau)}{d \tau^n}d\tau ,~~ t>a,\,\,
\alpha \in [n-1,n).
\end{equation}
A third alternative definition is the Gr\"{u}nwald-Letnikov
derivative,  based on finite differences to generalize the classical
derivative. This
 is equivalent to the Riemann-Liouville derivative if one ignores assumptions on the regularity of the functions.

 For (\ref{RLDeri}), the following holds \cite{Deng:07}
 $$\lim\limits_{\alpha \rightarrow
(n-1)^+}{_aD_t^\alpha}x(t)=\frac{d^{n-1}x(t)}{dt^{n-1}},~~
\lim\limits_{\alpha \rightarrow
n^-}{_aD_t^\alpha}x(t)=\frac{d^nx(t)}{dt^n};$$ and for
(\ref{DefFracInt3}), $\lim\limits_{\alpha \rightarrow
n}{_aD_t^{-\alpha}}x(t)={_aD_t^{-n}}x(t)$. Hence, the operator
${_aD_t^\alpha} \,(\alpha \in {R})$ makes sense and `connects' the
orders of the fractional derivative from $-\infty$ to $+\infty$.
However, for (\ref{CDeri}), one recovers
$$\lim\limits_{\alpha \rightarrow
(n-1)^+}{_a^CD_t^\alpha}x(t)=\frac{d^{n-1}x(t)}{dt^{n-1}}-\frac{d^{n-1}x(t)}{dt^{n-1}}\big|_{t=a},~~
\lim\limits_{\alpha \rightarrow
n^-}{_a^CD_t^\alpha}x(t)=\frac{d^nx(t)}{dt^n}.$$ The
Riemann-Liouville derivative and the Caputo derivative are
equivalent if $x(t)$ is sufficiently smooth and satisfies
$x^{(j)}(a)=0,\,j=0,2,\cdots,n-1.$

In this work, we develop a  local discontinuous Galerkin (LDG)
methods for the one-term and multi-term initial value problems for
fractional ordinary differential equations (FODEs). For ease of
presentation, we shall focus the discussion on the following two
types of FODEs
\begin{equation}\label{FirstClass}
{_a^CD_t^\alpha}x(t)=f(x,t),
\end{equation}
and
\begin{equation}\label{SecondClass}
{_a^CD_t^\alpha}x(t)+d(t) \frac{d^mx(t)}{dt^m}=f(x,t),
\end{equation}
where $\alpha \in (0,1] $, and $m$ is a positive integer and also
include examples of $\alpha \in [1,2]$ in the interest of
generalization. For (\ref{FirstClass}) and (\ref{SecondClass}), the
initial conditions can be specified exactly as for the classical
ODEs, i.e., the values of $x^{(j)}(a)$ must be given, where
$j=0,1,2,\cdots,\lfloor \alpha\rfloor$ for (\ref{FirstClass}) and
$j=0,1,2,\cdots,\max\{\lfloor \alpha\rfloor,m-1\}$ for
(\ref{SecondClass}).

It appears that the earliest numerical methods used in the
engineering community for such problems are the predictor-corrector
approach originally presented in \cite{Diethelm:02}, later slightly
improved in \cite{Deng:07}, and a method using a series of classical
derivatives to approximate the fractional derivative, realized by
using frequency domain techniques based on Bode diagrams
\cite{Hartley:08}. For the second method, ways to evaluate  the time
domain error introduced in the frequency domain approximations
remains open.

The discontinuous Galerkin (DG) methods have been well developed to
solve classical differential equations \cite{Hesthaven:08},
initiated for the classical ODEs \cite{Delfour:81} with substantial
later work, mostly related to discontinuous Galerkin methods for the
related Volterra integro-differential equation, including a priori
analysis \cite{Schoetzau:00}, $hp$-adaptive methods
\cite{Brunner:06,Mustapha:11} and recent work on super convergence
in the $h$-version \cite{Mustapha:13}. This has been extended to
approximate the fractional  spatial derivatives \cite{Deng:11} to
solve fractional diffusion equation by using the idea of local
discontinuous Galerkin (LDG) methods \cite{Bassi:97,Cockburn:98}. In
this work, we discuss DG methods to allow for the approximate
solution of general FODEs. All the advantages/characteristics of the
spatial DG methods  carries over to this case with a central one
being the ability  to solve the equation interval by interval when
the upwind flux, taking the value of $x(t)$ at a discontinuity point
$t_j$ as $x(t_j^-)$, is used. However, this is a natural choice,
since for the initial value problems for the fractional (or
classical) ODEs, the information travels forward in time. This
implies that we just invert a local low order matrix rather than a
global {\em full} matrix.  The LDG methods for the first order FODEs
(\ref{FirstClass}) have optimal order of convergence $k+1$ in the
$L^2$ norm and we observe superconvergence of order
$k+1+\min\{k,\alpha\}$ at the downwind point of each element. Here
$k$ is the degree of the approximation polynomial used in an
element. For the two-term FODEs (\ref{SecondClass}), the LDG methods
retains optimal convergence order in $L^2$-norm, and
superconvergence at the downwind point of each element as
$k+1+\min\{k,\max\{\alpha,m\}\}$. We shall discuss the underlying
mechanism of this  superconvergence and illustrate the results of
the analysis through a number of examples, including some going
beyond the theoretical developments presented here.

What remains of the paper is organized as follows. In Section 2, we
present the LDG schemes for the FODEs, discuss the numerical
stability of the scheme for the linear case of (\ref{FirstClass})
with $\alpha \in [0,1]$, and uncover the mechanism of
superconvergence. The analysis is supported by a number of
computational experiments and we also illustrate how the proposed
scheme can be applied to compute the generalized Mittag-Leffler
functions. In Section 3 we discuss a number of generalizations of
the scheme, illustrated by a selection of numerical examples and
Section 4 contains a few concluding remarks.

\section{LDG schemes for the FODEs}
The basic idea in the design the LDG schemes is to rewrite the FODEs
as a system of the first order classical ODEs and a fractional
integral. Since the integral operator naturally connects the
discontinuous function, we need not add a penalty term or introduce
a numerical fluxes for the integral equation. However,  for the
first order ODEs, upwind fluxes are used. In this section, we
present the LDG schemes, prove numerical stability and discuss the
underlying mechanism of superconvergence.

\subsection{LDG schemes}
We consider (\ref{FirstClass}) and rewrite it as
\begin{equation} \label{ModelILDG}
\left\{\begin{array}{l} \displaystyle x_1(t)-\frac{dx_0(t)}{dt}=0
,~~~ {_0D_t^{-(1-\alpha)}}x_1(t)=f(x,t),~~~ \alpha \in (0,1],
\\
\\
x_0(0)=x_0.
\end{array}
\right.
\end{equation}
We consider a scheme for solving (\ref{ModelILDG}) in the interval
$\Omega=[0,T]$. Given the nodes $0=t_0<t_1<\cdots<t_{n-1}<t_n=T$, we
define the mesh $\mathcal{T}=\{I_j=(t_{j-1},t_j),\,j=1,2,\cdots,n\}$
and set $h_j:=|I_j|=t_j-t_{j-1}$ and $h:=\max_{j=1}^n h_j$.
Associated with the mesh $\mathcal{T}$, we define the broken Sobolev
spaces
$$
L^2(\Omega,\mathcal{T}):=\{ v: \Omega \rightarrow {R}\, \big| \,
v|_{I_j} \in L^2(I_j),\,j=1,2,\cdots,n \};
$$
and
$$
H^1(\Omega,\mathcal{T}):=\{ v: \Omega \rightarrow {R}\, \big| \,
v|_{I_j} \in H^1(I_j),\,j=1,2,\cdots,n \}.
$$
For a function $v \in H^1(\Omega,\mathcal{T})$, we denote the
one-sided limits at the nodes $t_j$ by
$$
v^{\pm}(t_j)=v(t_j^{\pm}):=\lim\limits_{t \rightarrow t_j^\pm}v(t).
$$
Assume that the solutions belong to the corresponding spaces:
$$(x_0(t),x_1(t)) \in H^1(\Omega,\mathcal{T}) \times
L^2(\Omega,\mathcal{T}).$$ We further define $X_i$ as the
approximation functions of $x_i$ respectively, in the finite
dimensional subspace $V \subset H^1(\Omega, \mathcal{T})$; and
choose $V$ to be the space of discontinuous, piecewise polynomial
functions
$$ V=\{v:\, \Omega \rightarrow {R}\, \big|\,v|_{I_j} \in
\mathcal{P}^k(I_j),\, j=1,2,\cdots,n \},
$$
where $\mathcal{P}^k(I_j)$ denotes the set of all polynomials of
degree less than or equal to $k$ on $I_j$. Using the upwind fluxes
for the first order classical ODEs and discretizing the integral
equation,  we seek $X_i \in V$ such that for all $v_i \in V$, and
$j=1,2,\cdots, n$, the following holds

\begin{equation}\label{Scheme1}
\left\{\begin{array}{l}
\displaystyle\big(X_1,v_0\big)_{I_j}+\bigg(X_0,\frac{dv_0}{dt}\bigg)_{I_j}-\big(
X_0(t_j^-)v_0(t_j^-)-X_0(t_{j-1}^-)v_0(t_{j-1}^+)\big)=0,
\\
\\
\displaystyle \big( {_0D_t^{-(1-\alpha)}}X_1,v_1
\big)_{I_j}=\big(f(X_0,t),v_1\big)_{I_j},
\\
\\
X_0(t_0^-)=x_0.
\end{array}
\right.
\end{equation}

{\em Remark:} If take $\alpha \in (0,1]$ and $m=1$ in
(\ref{SecondClass}), the scheme is
$$
\left\{\begin{array}{l}
\displaystyle\big(X_1,v_0\big)_{I_j}+\bigg(X_0,\frac{dv_0}{dt}\bigg)_{I_j}-\big(
X_0(t_j^-)v_0(t_j^-)-X_0(t_{j-1}^-)v_0(t_{j-1}^+)\big)=0,
\\
\\
\displaystyle \big( {_0D_t^{-(1-\alpha)}}X_1+X_1,v_1
\big)_{I_j}=\big(f(X_0,t),v_0\big)_{I_j},
\\
\\
X_0(t_0^-)=x_0.
\end{array}
\right.
$$
The scheme is clearly consistent, i.e,, the exact solutions of the
corresponding models satisfy (\ref{Scheme1}). Furthermore, since an
upwind flux is used  the solutions can be computed interval by
interval and if $f(x,t)$ is a linear function in $x$, we just need
to invert a small matrix in each interval to recover the solution.

\subsection{Numerical stability}
We consider the question of stability for the linear case of
(\ref{Scheme1}):

\begin{equation}
\label{LinearModelI} \left\{\begin{array}{l} \displaystyle
{_0^CD_t^\alpha}x(t)=Ax(t)+B(t),~~~~~~ \alpha \in (0,1),
\\
\\
x(0)=x_0,
\end{array}
\right.
\end{equation}
where $A$ is a negative constant and $B(t)$ is sufficiently regular
to ensure existence and uniqueness. The numerical scheme of
(\ref{LinearModelI}) is to find $X_i \in V$ such that
\begin{equation} \label{SchemeLinearModelI}
\left\{
\begin{array}{l}
\displaystyle \big(X_0(t_j^-)v_0(t_j^-)-X_0(t_{j-1}^-)v_0(t_{j-1}^+)
\big)-\big(X_1,v_0 \big)_{I_j}
-\bigg(X_0,\frac{dv_0}{dt}\bigg)_{I_j}=0,
\\
\\
\big( {_0D_t^{-(1-\alpha)}}X_1,v_1
\big)_{I_j}=\big(AX_0+B,v_1\big)_{I_j},
\\
\\
X_0(t_0^-)=x_0,
\end{array}
\right.
\end{equation}
holds for all $v_i \in V$. First we present a lemma here. Based on
the semigroup property of fractional integral operators and Lemmas
2.1 and 2.6 in \cite{Deng:11}, we have

\begin{lemma}\label{LemmaNorm} For $\beta>0$, $\alpha \in (0,1)$,
\begin{equation}\label{Commute}
\big({_0D_t^{-\beta}}u,v
\big)_{L^2([0,t_j])}=\big(u,{_tD_{t_j}^{-\beta}}
v\big)_{L^2([0,t_j])};
\end{equation}
\begin{equation}\label{LemmaNormEq}
\begin{array}{ll}
\displaystyle\big({_0D_t^{\alpha-1}}v, v\big)_{L^2([0,t_j])} &
=\big({_0D_t^{\frac{\alpha-1}{2}}}v,
{_tD_{t_j}^{\frac{\alpha-1}{2}}}v\big)_{L^2([0,t_j])}
\\
\\
&\displaystyle=\cos\bigg( \frac{(\alpha-1)\pi}{2}
\bigg)\|v\|_{H^{\frac{\alpha-1}{2}}([0,t_j])}^2.
\end{array}
\end{equation}
\end{lemma}
Let ${\widetilde X}_i \in V $ be the approximate solution of $X_i$
and denote $e_{X_i}:={\widetilde X}_i-X_i$ as the numerical errors.
Stability of (\ref{SchemeLinearModelI}) is established in the
following theorem.
\begin{theorem}[$L^\infty$ stability] \label{StabilityTheorem}
Scheme (\ref{SchemeLinearModelI}) is $L^\infty$ stable; and the
numerical errors satisfy
\begin{equation}\label{Stability}
~~e_{X_0}^2(t_n^-)=e_{X_0}^2(t_0^-)-\sum\limits_{i=1}^n(\delta
e_{X_0}(t_{i-1}))^2+\frac{2}{A}\cos\bigg( \frac{(\alpha-1)\pi}{2}
\bigg)\|e_{X_1}\|_{H^{\frac{\alpha-1}{2}}([0,t_j])}^2,
\end{equation}
where $\delta
e_{X_i}(t_{i-1})=e_{X_i}(t_{i-1}^+)-e_{X_i}(t_{i-1}^-)$.
\end{theorem}
Since $A<0$, the scheme is dissipative.

\begin{proof}
From (\ref{SchemeLinearModelI}), we recover the error equation
\begin{equation} \label{Perturbation}
\left\{
\begin{array}{l}
\displaystyle
\big(e_{X_0}(t_i^-)v_0(t_i^-)-e_{X_0}(t_{i-1}^-)v_0(t_{i-1}^+)
\big)-\big(e_{X_1},v_0 \big)_{I_i}
-\bigg(e_{X_0},\frac{dv_0}{dt}\bigg)_{I_i}=0,
\\
\\
-\frac{1}{A}\big( {_0D_t^{-(1-\alpha)}}e_{X_1},v_1
\big)_{I_i}=-\big(e_{X_0},v_1\big)_{I_i},
\end{array}
\right.
\end{equation}
for all $v_i \in V$. Taking $v_0=e_{X_0}$, $v_1=e_{X_1}$, and adding
the two equations, we obtain
$$
\big( e_{X_0}^2(t_i^-) -e_{X_0}(t_{i-1}^-)e_{X_0}(t_{i-1}^+)
\big)-\frac{1}{2}\big(e_{X_0}^2(t_i^-) -e_{X_0}^2(t_{i-1}^+)
\big)-\frac{1}{A}\big({_0D_t^{-(1-\alpha)}}e_{X_1}
,e_{X_1}\big)_{I_i}=0.
$$
Summing the equations for $i=1,2,\cdots,n$ leads to
$$
\begin{array}{l}
\displaystyle
e_{X_0}^2(t_n^-)-e_{X_0}^2(t_0^-)+\sum\limits_{i=1}^n\big(e_{X_0}^2(t_{i-1}^-)-2e_{X_0}(t_{i-1}^-)e_{X_0}(t_{i-1}^+)+e_{X_0}^2(t_{i-1}^+)
\big)
\\
\\
\displaystyle-\frac{2}{A}
\big({_0D_t^{-(1-\alpha)}}e_{X_1},e_{X_1}\big)_{[0,I_n]}=0,
\end{array}
$$
and
$$
e_{X_0}^2(t_n^-)=e_{X_0}^2(t_0^-)-\sum\limits_{i=1}^n\big(e_{X_0}(t_{i-1}^+)-e_{X_0}(t_{i-1}^-)
\big)^2+\frac{2}{A}
\big({_0D_t^{-(1-\alpha)}}e_{X_1},e_{X_1}\big)_{[0,I_n]}.
$$
Using (\ref{LemmaNormEq}) of Lemma \ref{LemmaNorm} yields the
desired result.
\end{proof}

\subsection{Mechanism of superconvergence}
As we shall see shortly, the  proposed LDG scheme is $k+1$ optimally
convergent in the $L^2$-norm but superconvergent at downwind points.
This is also known for the classic case where the downwind
convergence is $2k+1$ \cite{Delfour:81,Adjerid:02}. However, for the
fractional case,  the order of the superconvergence depends on the
order of the Caputo fractional derivatives. To understand this, let
us again focus on the linear case of (\ref{Scheme1}).

The error equation corresponding to (\ref{SchemeLinearModelI}) is
\begin{equation} \label{Error}
\left\{
\begin{array}{l}
\displaystyle
\big(E_{X_0}(t_i^-)v_0(t_i^-)-E_{X_0}(t_{i-1}^-)v_0(t_{i-1}^+)
\big)-\big(E_{X_1},v_0 \big)_{I_i}
-\bigg(E_{X_0},\frac{dv_0}{dt}\bigg)_{I_i}=0,
\\
\\
-\frac{1}{A}\big( {_0D_t^{-(1-\alpha)}}E_{X_1},v_1
\big)_{I_i}=-\big(E_{X_0},v_1\big)_{I_i},
\end{array}
\right.
\end{equation}
where $E_{X_i}=x_i-X_i$. In (\ref{Error}), taking $v_i$ to be
continuous on the interval $[0,t_j]$ and summing the equations for
$i=1,2,\cdots,n$ lead to
\begin{equation}\label{supercov1}
\begin{split}
&E_{X_0}(t_n^-)v_0(t_n^-)-(E_{X_1},v_0)_{[0,t_n]}-\left(E_{X_0},\frac{dv_0}{dt}\right)_{[0,t_n]}\\
&\quad-\frac{1}{A}({_0D_t^{-(1-\alpha)}}E_{X_1},v_1)_{[0,t_n]}+(E_{X_0},v_1)_{[0,t_n]}=0.
\end{split}
\end{equation}
Following (\ref{Commute}), we rewrite this as
\begin{equation}\label{supercov2}
\begin{split}
&E_{X_0}(t_n^-)v_0(t_n^-)-(E_{X_1},v_0)_{[0,t_n]}-\left(E_{X_0},\frac{dv_0}{dt}\right)_{[0,t_n]}\\
&\quad-\frac{1}{A}(E_{X_1},{_tD_{t_n}^{-(1-\alpha)}}v_1)_{[0,t_n]}+(E_{X_0},v_1)_{[0,t_n]}=0.
\end{split}
\end{equation}
Rearranging the terms of (\ref{supercov2}) results in
\begin{equation}\label{supercov3}
E_{X_0}(t_n^-)v_0(t_n^-)-\left(E_{X_1},v_0+\frac{1}{A}{_tD_{t_n}^{-(1-\alpha)}}v_1
\right)_{[0,t_n]}-\left(E_{X_0}, \frac{dv_0}{dt}-v_1
\right)_{[0,t_n]}=0.
\end{equation}
Solving
$\tilde{v}_0+\frac{1}{A}{_tD_{t_n}^{-(1-\alpha)}}\tilde{v}_1=0$ and
$\frac{d\tilde{v}_0}{dt}-\tilde{v}_1=0$ for $t \in [0,t_n]$ with
$\tilde{v}_0(t_n)=E_{X_0}(t_n^-)$, we get
\begin{equation}\label{supercov4}
\tilde{v}_0(t)=(E_{X_0}(t_n^-)/E_{\alpha,1}(-At_n^\alpha))
E_{\alpha,1}(-At^\alpha),
\end{equation}
where $E_{\alpha,1}$ is the Mittag-Leffler function. Taking $v_i$ as
the $L^2$ projection of $\tilde{v}_i$ onto $\mathcal{P}^k$, we
recover that if $\alpha$ is an integer, there  exists
\begin{equation}\label{supercov5}
\left\|v_0+\frac{1}{A}{_tD_{t_n}^{-(1-\alpha)}}v_1+\frac{dv_0}{dt}-v_1\right\|_{L^2([0,t_n])}=O(h^k);
\end{equation}
due to the regularity of the Mittag-Leffler function
\cite{Mainardi:14}. For the fractional case, the approximation
\cite{Guo:06} yields
\begin{equation}\label{supercov6}
\left\|v_0+\frac{1}{A}{_tD_{t_n}^{-(1-\alpha)}}v_1+\frac{dv_0}{dt}-v_1\right\|_{L^2([0,t_n])}=O(h^{\min\{k,\alpha\}}).
\end{equation}
Combined with classic polynomial approximation results for $E_{X_i}$
\cite{Delfour:81}, yields an order of convergence at the downwind
point of $k+1+\min \{ k, \alpha \}$.

\subsection{Numerical experiments}
Let us consider a few numerical examples to qualify the above
analysis.

All results assume that the corresponding analytical solutions are
sufficiently regular. For showing the effectiveness of the LDG
schemes and further confirming the predicted convergence orders,
both linear and nonlinear cases are considered. Finally, we shall
also consider the computation of the generalized Mittag-Leffler
functions using the LDG scheme.

We use Newton's method to solve the nonlinear systems. The initial
guess in the interval $I_j$ ($j \geq 2$) is given as
$X_i(t_{j-1}^-)$ or by extrapolating forward to the interval $I_j$.
For the interval $I_1$, we use $X_0(t_0^-)$ $(=x_0)$ as initial
guess.

\subsubsection{Numerical results for (\ref{FirstClass}) with $\alpha \in (0,1]$}

We first consider examples to confirm that the convergence order of
(\ref{Scheme1}) is $k+1+\alpha$ at downwind points and $k+1$ in the
$L^2$ sense, respectively, where $k$ is the degree of the polynomial
used in an element and $\alpha \in [0,1)$. However, when $\alpha=1$
the convergence order is $2k+1$ at downwind point and
still $k+1$ in $L^2$ sense in agreement with classic theory \cite{Delfour:81}. \\

{\em Example L1.} On the computational domain $t\in \Omega=(0,1)$,
we consider
\begin{equation} \label{ExampleL1}
{_0^CD_t^\alpha
x(t)}=-2x(t)+\frac{\Gamma(6)}{\Gamma(6-\alpha)}t^{5-\alpha}+2t^5+2,~~~
\alpha \in [0,1],
\end{equation}
with the initial condition $x(0)=1$ and the exact solution
$x(t)=t^5+1$. Note that when $\alpha=0$ this condition is still
required for the form of (\ref{ModelILDG}).

\begin{figure}
\vbox{ \hbox{
\includegraphics[width=2.5in,angle=0]{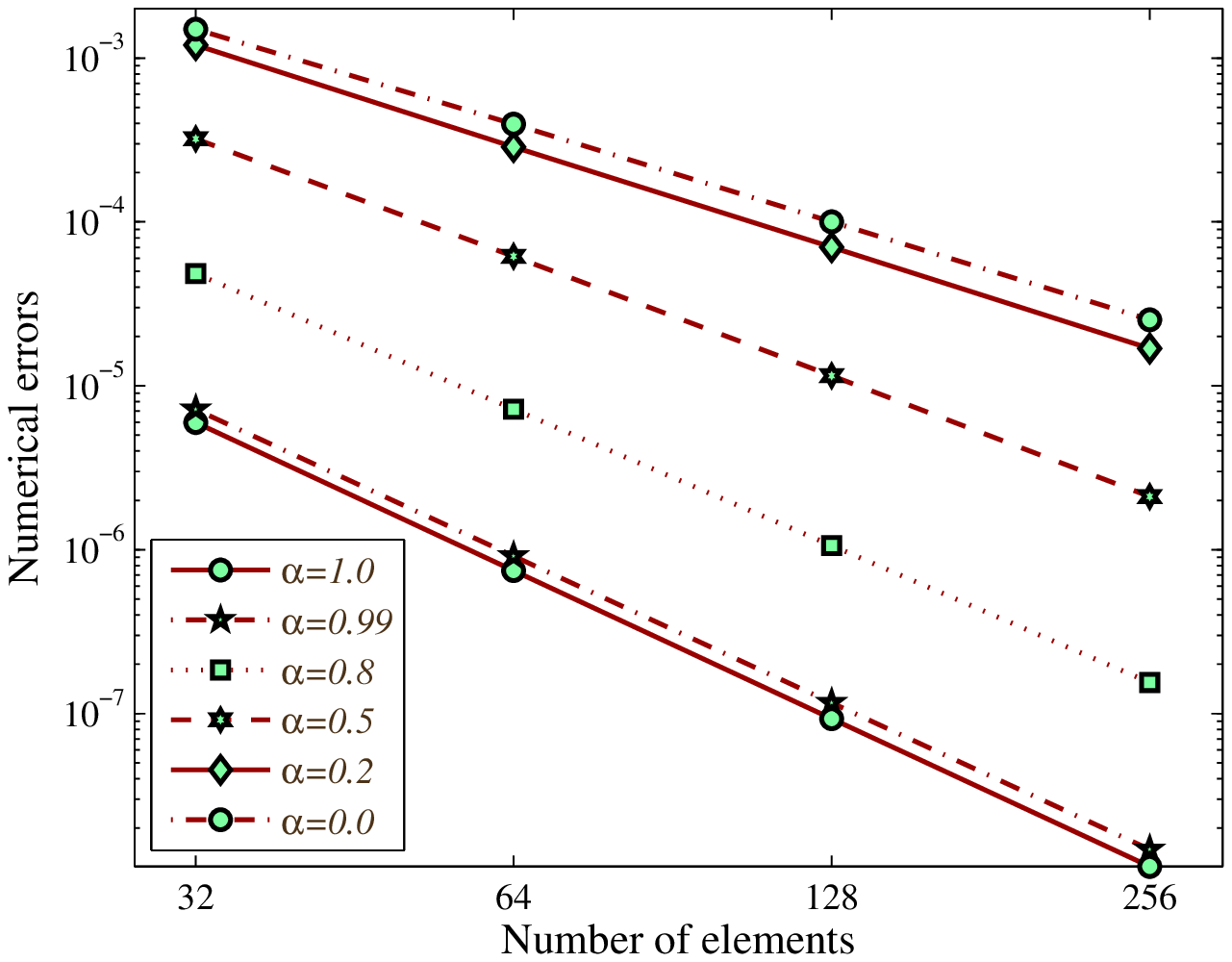}
\includegraphics[width=2.5in,angle=0]{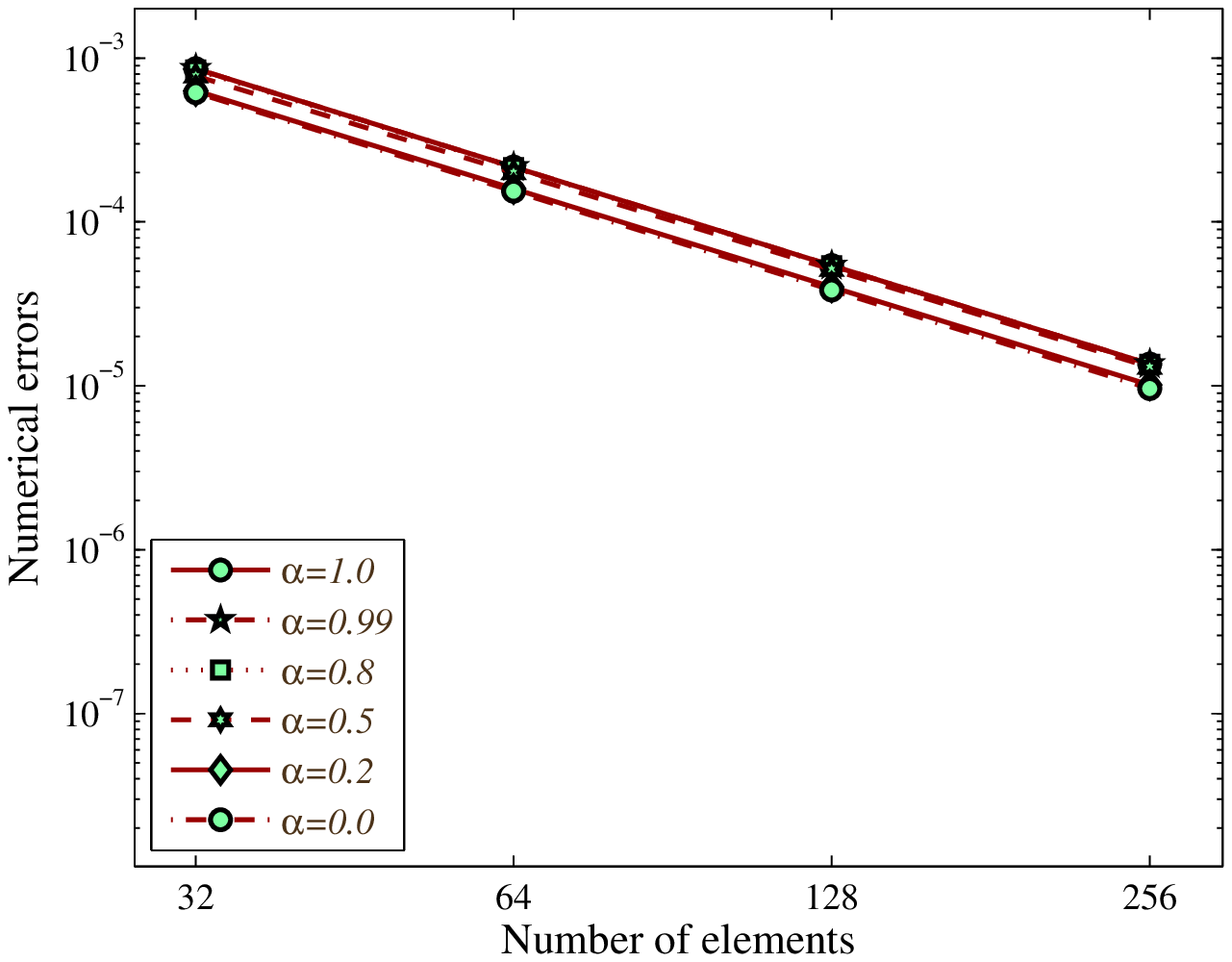}}
\hbox{
\includegraphics[width=2.5in,angle=0]{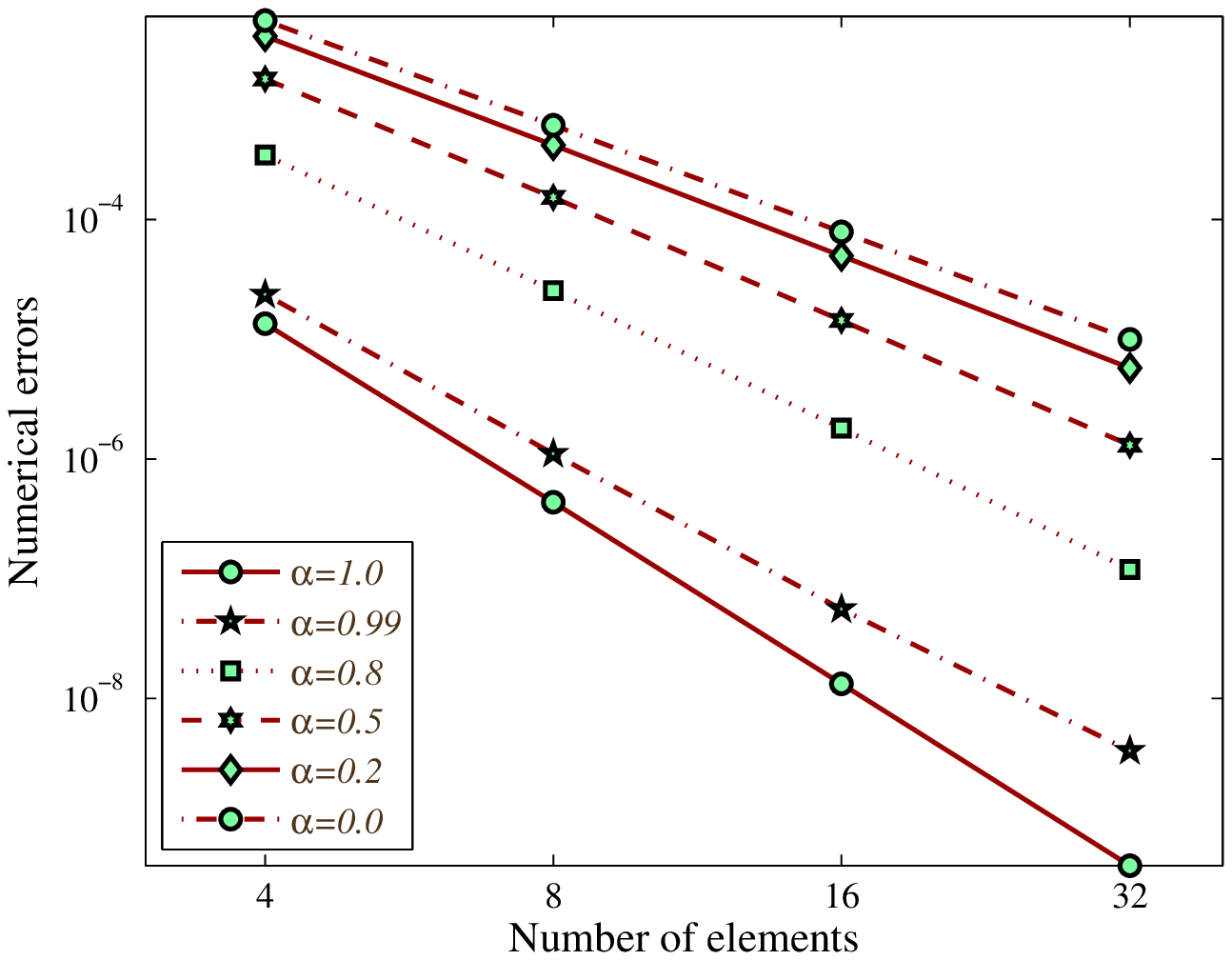}
\includegraphics[width=2.5in,angle=0]{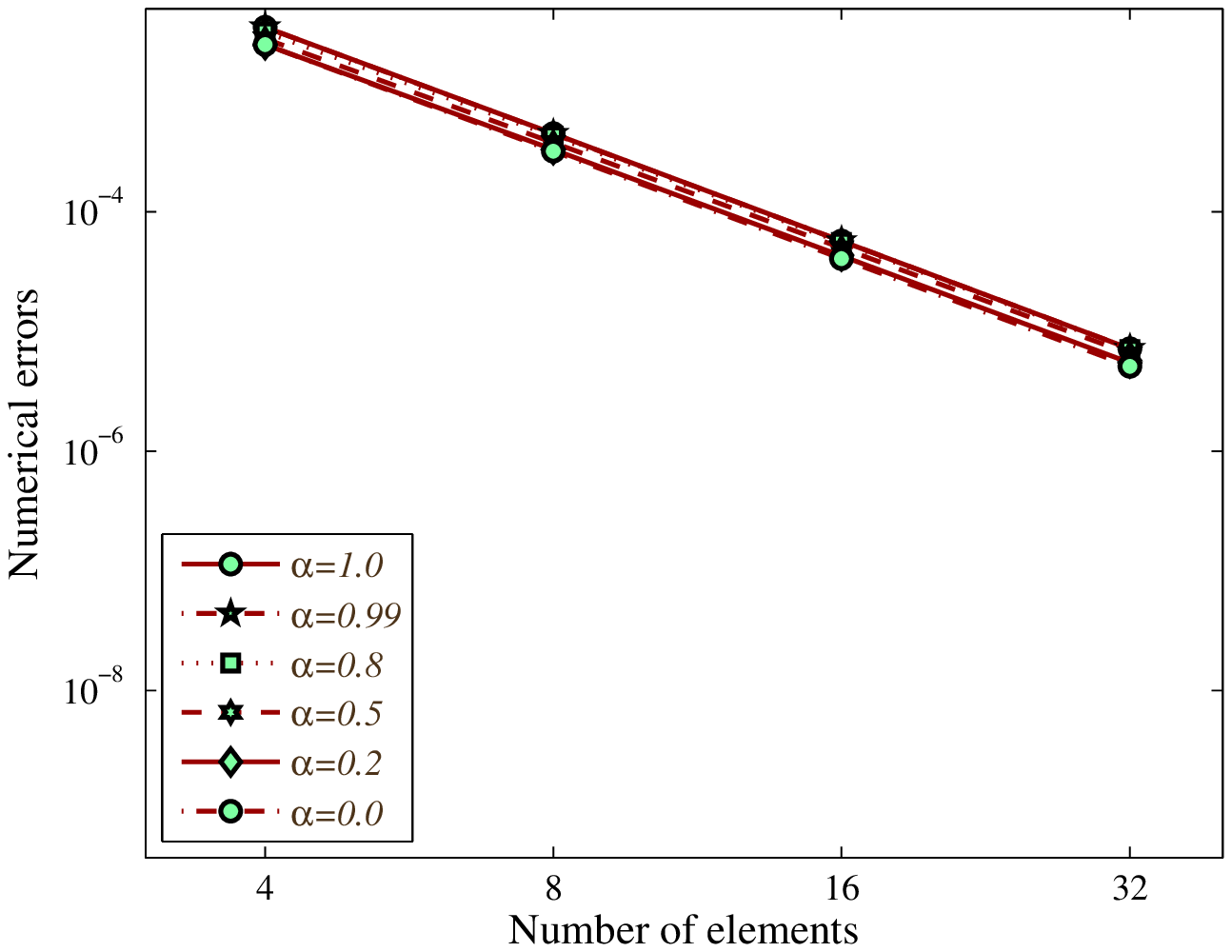}}
\hbox{
\includegraphics[width=2.5in,angle=0]{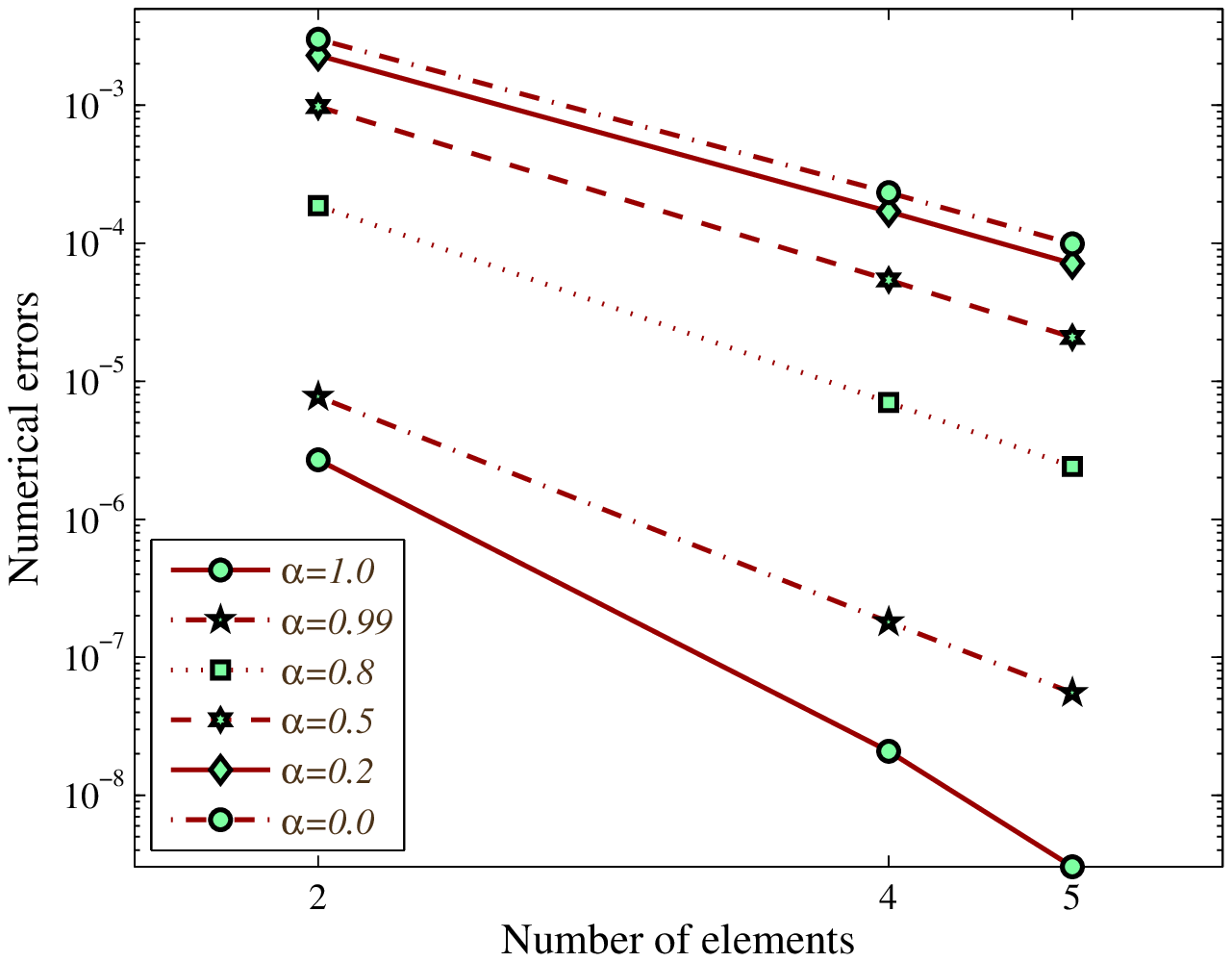}
\includegraphics[width=2.5in,angle=0]{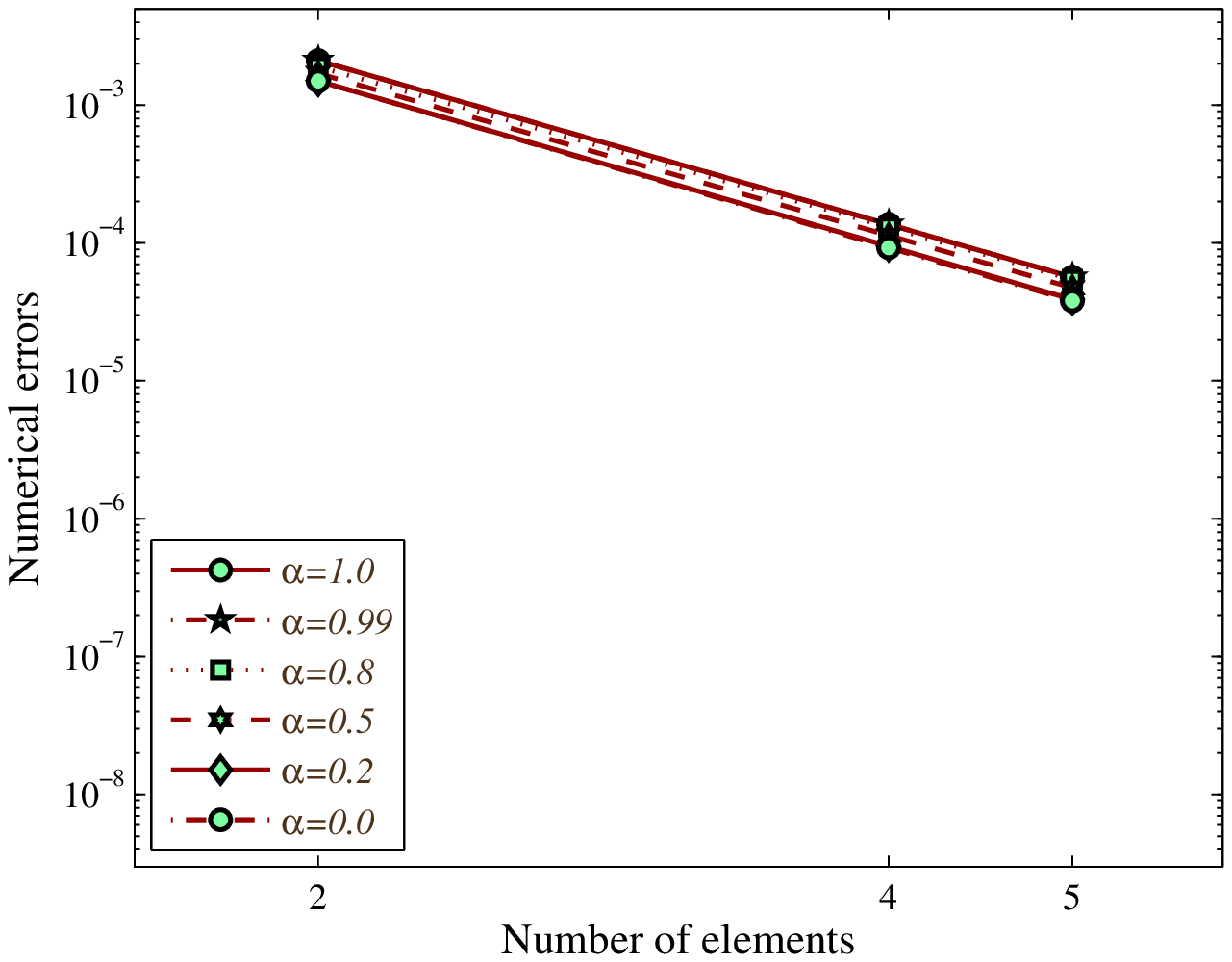}}}
\caption{The convergence of (\ref{Scheme1}) for $k=1$ (top row),
$k=2$ (middle row) and $k=3$ (bottom row) in (\ref{ExampleL1}). In
the left column we show the convergence at the downwind points while
the right column displays $L^2$-convergence. \label{L1}}
\end{figure}
Figure \ref{L1} displays convergence the downwind point as well as
in $L^2$ for $k=1-3$, confirming optimal $L^2$ convergence and an
order of convergence
of $k+1+\alpha$ at the downwind point as predicted. \\
{\em Example N1.} We consider the nonlinear FODE on the domain $t
\in \Omega=(0,0.5)$,
\begin{equation}\label{ExampleN1}
{_0^CD_t^\alpha
x(t)}=-2x^2(t)+\frac{\Gamma(6)}{\Gamma(6-\alpha)}t^{5-\alpha}+2t^{10}+4t^5+2,~~~
\alpha \in [0,1],
\end{equation}
with the initial condition $x(0)=1$, and the exact solution
$x(t)=t^5+1$.
\begin{figure}
\vbox{ \hbox{
\includegraphics[width=2.5in,angle=0]{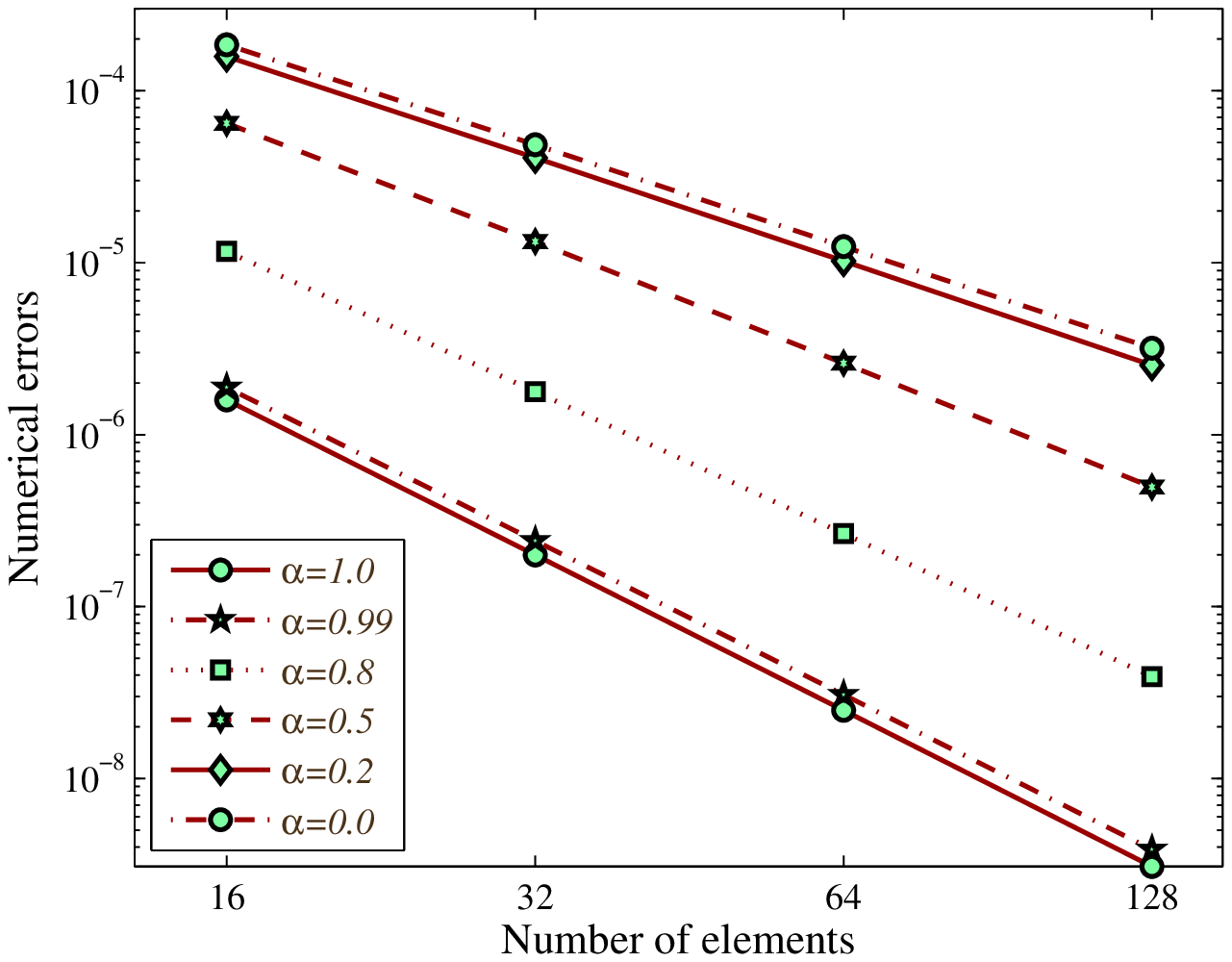}
\includegraphics[width=2.5in,angle=0]{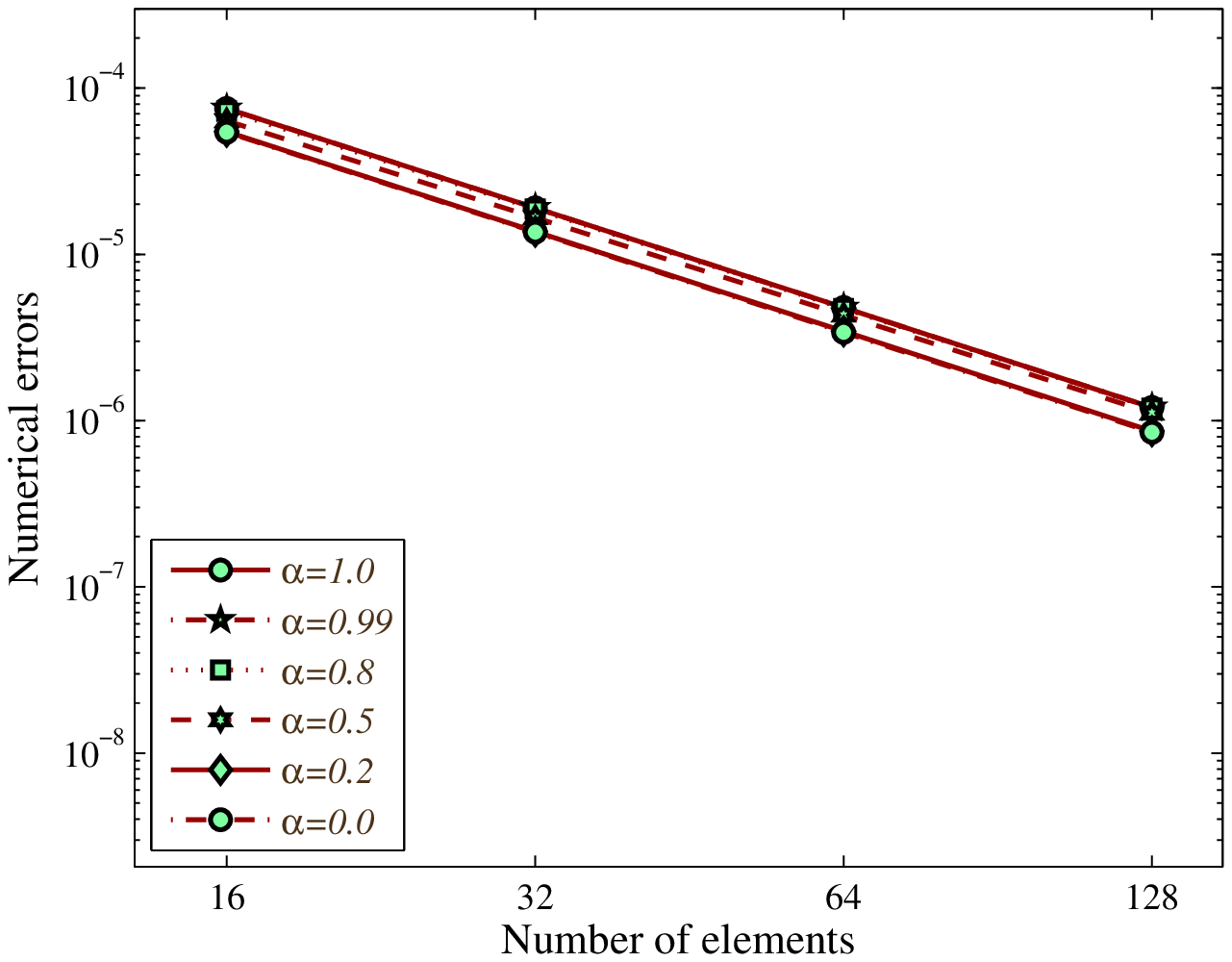}}
\hbox{
\includegraphics[width=2.5in,angle=0]{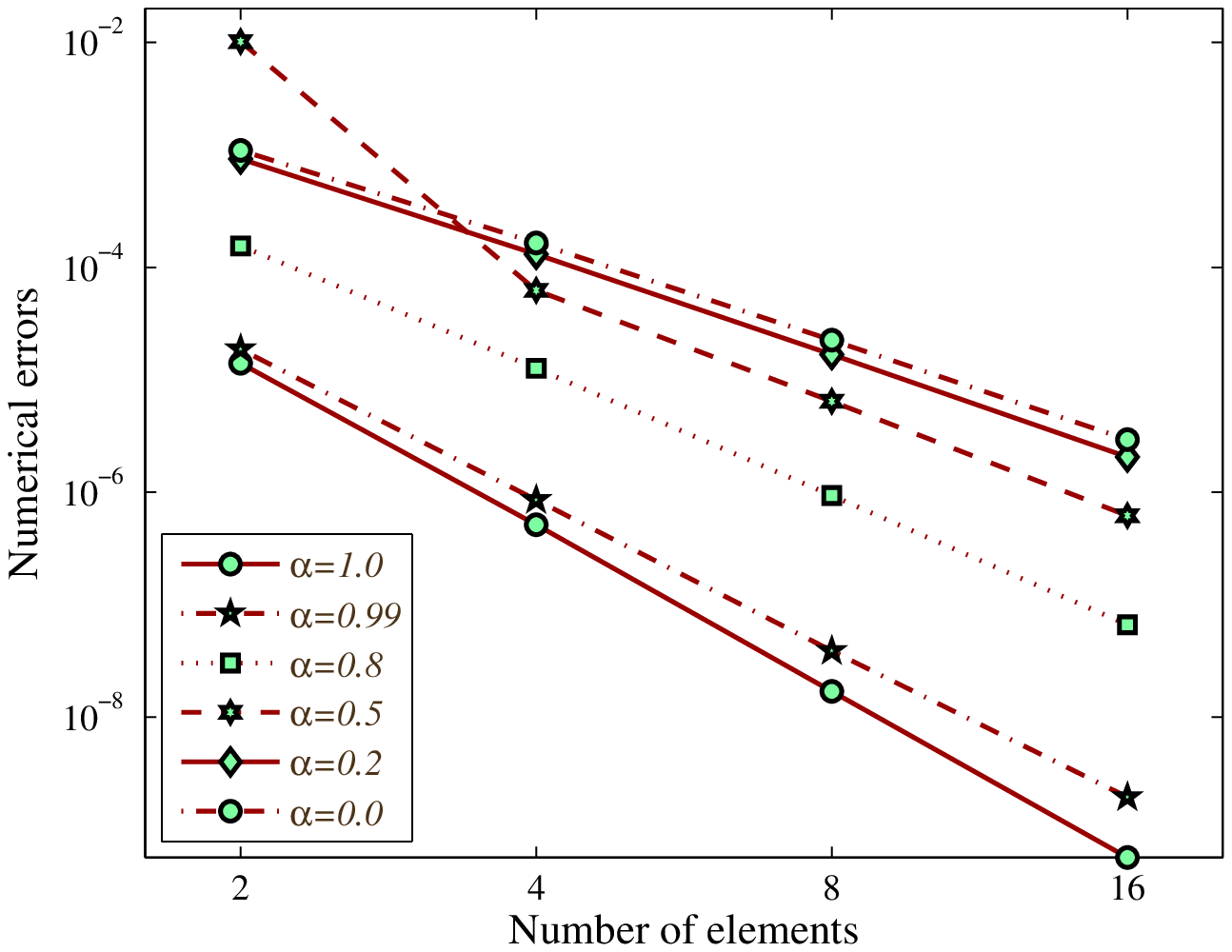}
\includegraphics[width=2.5in,angle=0]{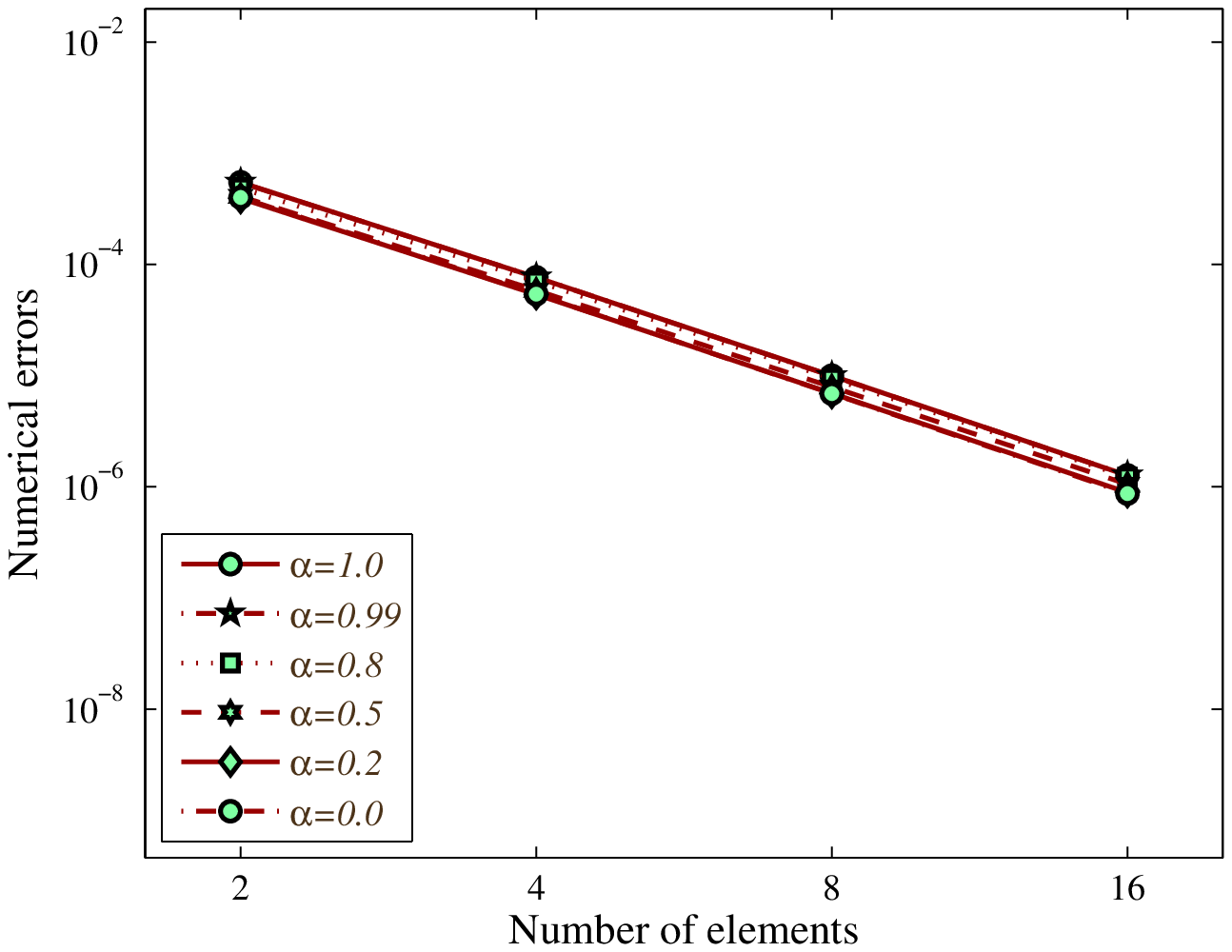}}
\hbox{
\includegraphics[width=2.5in,angle=0]{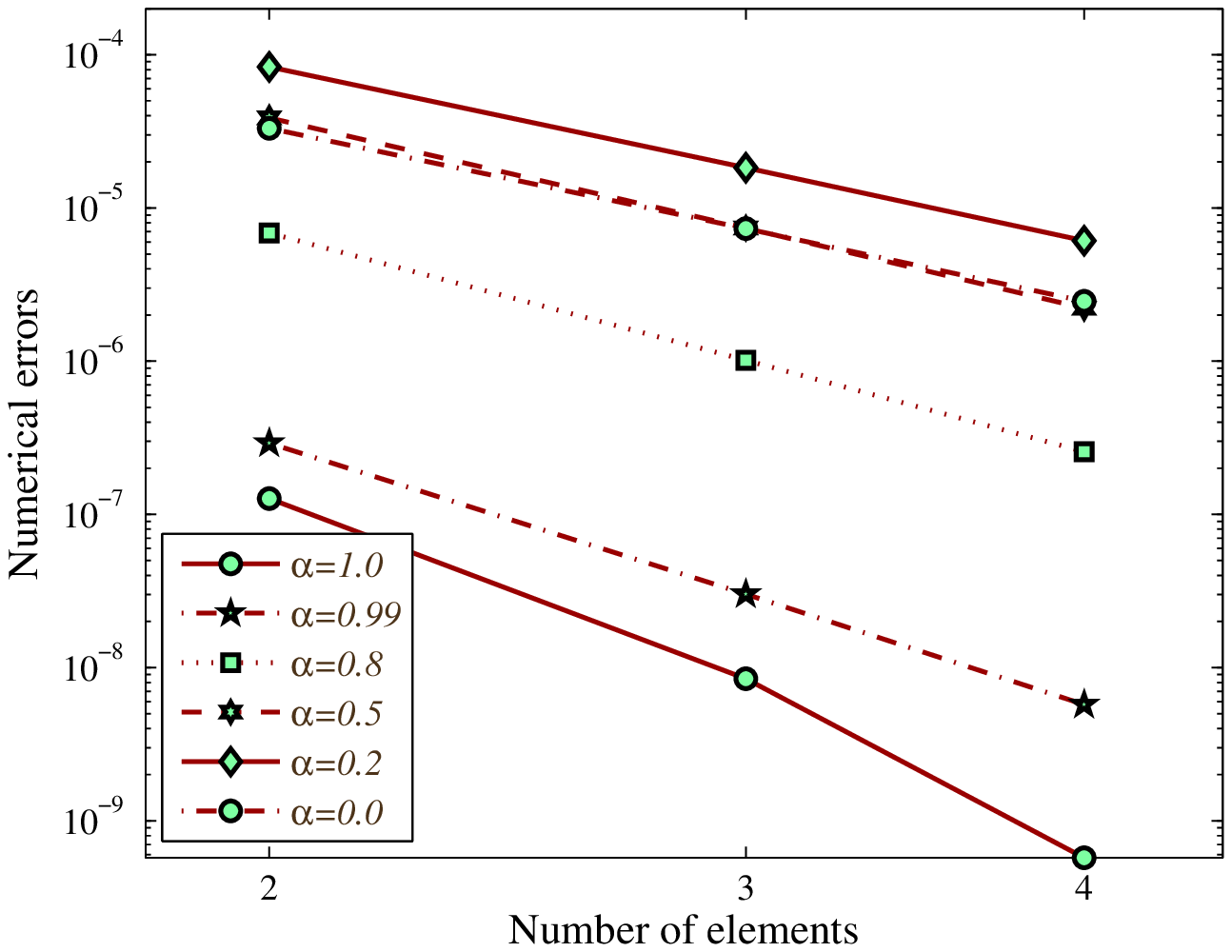}
\includegraphics[width=2.5in,angle=0]{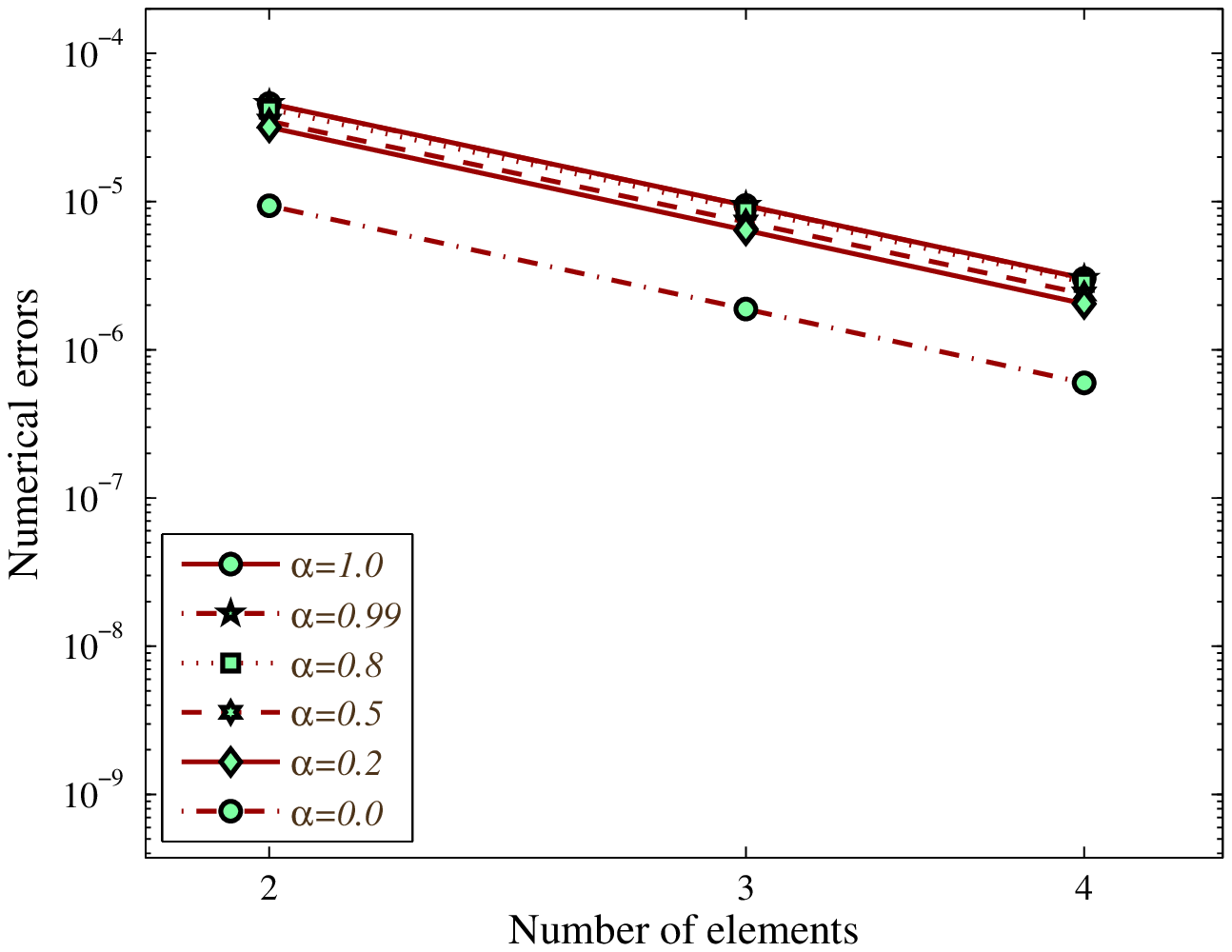}}}
\caption{The convergence of (\ref{Scheme1}) for $k=1$ (top row),
$k=2$ (middle row) and $k=3$ (bottom row) in (\ref{ExampleN1}). In
the left column we show the convergence at the downwind points while
the right column displays $L^2$-convergence. \label{N1}}
\end{figure}
The results in Figure \ref{N1} confirm that the optimal $L^2$
convergence and an order of convergence of $k+1+\alpha$ at the
downwind point carries over to the nonlinear case.

\subsubsection{Calculating the generalized Mittag-Leffler function}
\begin{figure}
\centering \hbox{
\includegraphics[width=2.5in]{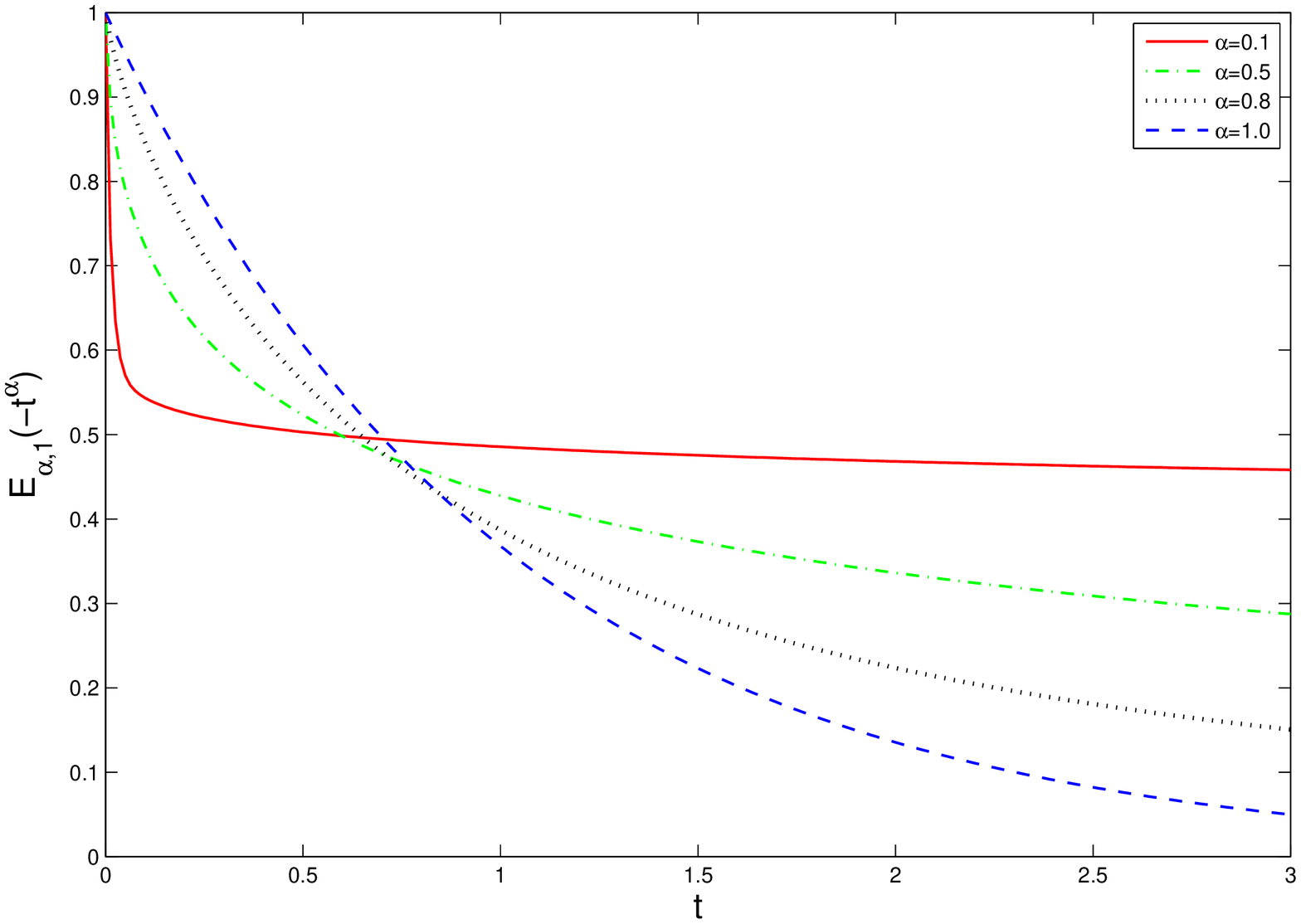}
\includegraphics[width=2.5in]{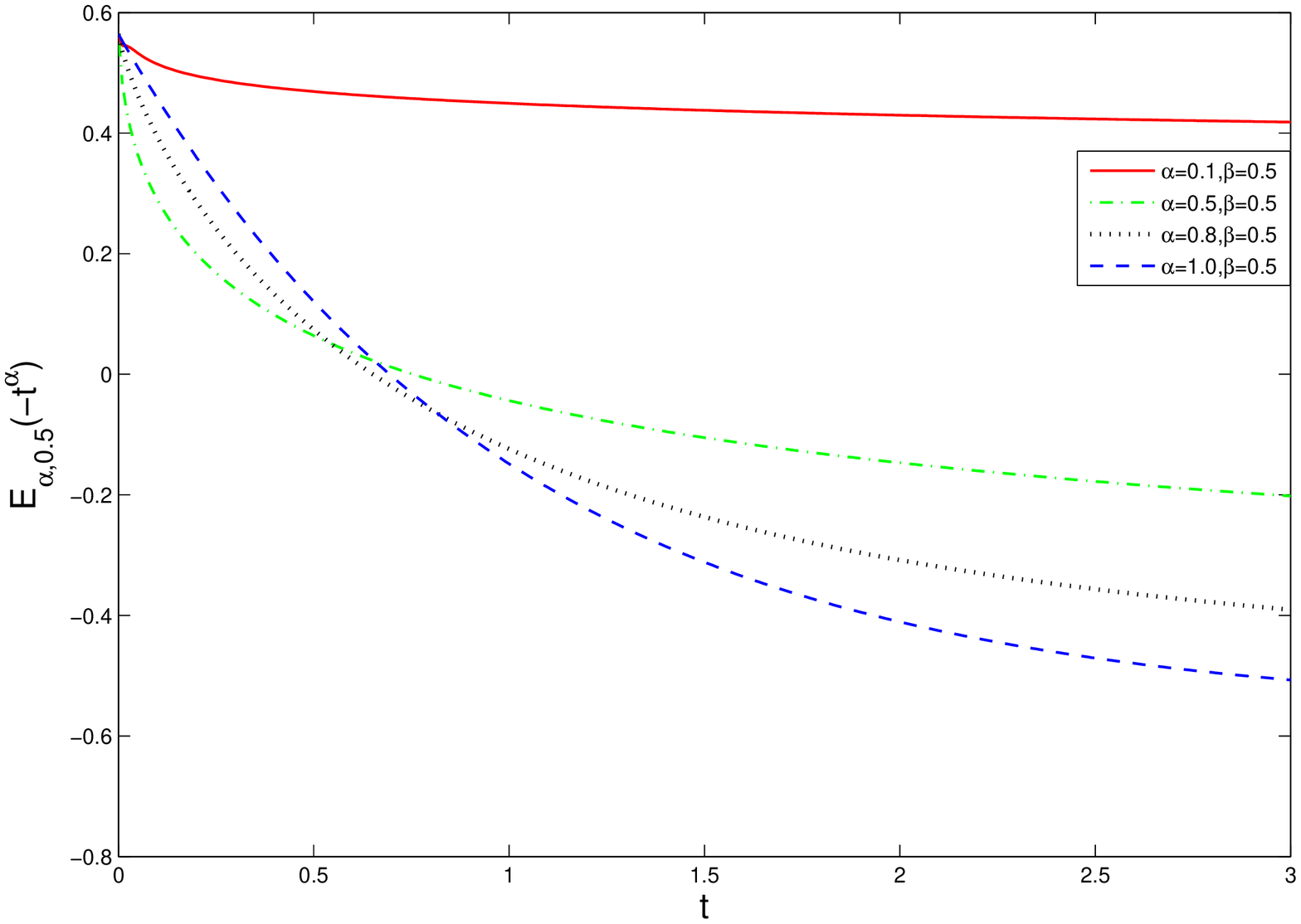}}
\caption{The Mittag-Leffler function $E_{\alpha,\beta}(-t^\alpha)$.
\label{fig}}
\end{figure}

As a more general example, let us use the efficient and accurate
solver for calculating the generalized Mittag-Leffler functions
defined as
\begin{equation}\label{MLF}
E_{\alpha,\beta}(At^\alpha)=\sum\limits_{k=0}^\infty
\frac{(At^\alpha)^k}{\Gamma(\alpha k+\beta)},~~ \Re(\alpha)>0.
\end{equation}
To build the relation between the Mittag-Leffler function and the
FODE consider
\begin{equation}\label{MLFEq1}
{_0^CD_t^\alpha}x(t)=Ax(t),~~
x(0)=1,\,x^\prime(0)=0,\cdots,x^{\lceil \alpha \rceil}(0)=0.
\end{equation}
Taking the Laplace transform on both sides of the above equation, we
recover
\begin{equation}\label{LT}
s^\alpha X(s)-s^{\alpha-1}x(0)=AX(s),
\end{equation}
where $X(s)$ is the Laplace transform of $x(t)$. From (\ref{LT}), we
 obtain
\begin{equation}\label{LTS}
X(s)=\frac{s^{\alpha-1}}{s^\alpha-A}.
\end{equation}
Using the inverse Laplace inverse transform in (\ref{LTS}) results
in
\begin{equation}\label{SEq}
x(t)=E_{\alpha,1}(At^\alpha)=\sum\limits_{k=0}^\infty
\frac{(At^\alpha)^k}{\Gamma(\alpha k+1)}.
\end{equation}
Since
$$
{_0D_t^{1-\beta}} \left( \frac{(At^\alpha)^k}{\Gamma(\alpha k+1)}
\right)=\frac{(At^\alpha)^k t^{\beta-1}}{\Gamma(\alpha k+\beta)},
$$
we recover
\begin{equation}\label{MLFEq2}
E_{\alpha,\beta}(At^\alpha)=t^{1-\beta}{_0D_t^{1-\beta}}
x(t)=t^{1-\beta}{_0D_t^{1-\beta}} \left(E_{\alpha,1}(At^\alpha)
\right).
\end{equation}
By solving (\ref{MLFEq1}) and  (\ref{MLFEq2}) or just (\ref{MLFEq1})
when $\beta=1$, we can efficiently calculate the generalized
Mittag-Leffler function $E_{\alpha,\beta}(At^\alpha)$ as illustrated
in Fig. \ref{fig}.

\section{Generalizations}
Let us now consider the generalized case, given as

\begin{equation}\label{SecondClass1}
{_a^CD_t^\alpha}x(t)+d(t) \frac{d^mx(t)}{dt^m}=f(x,t),
\end{equation}
where $\alpha $ in general is real, and $m$ is a positive integer.
We shall follow the same approach as previously, and consider the
system of equations

\begin{equation} \label{ModelVLDG}
\left\{\begin{array}{l}
\displaystyle x_{i+1}(t)-\frac{dx_i(t)}{dt}=0, i=0, \ldots, \max(m,p)-1, x_0(t) = x(t), \\
\\ \displaystyle  {_0D_t^{-(p-\alpha)}}x_p(t)+d(t) x_m(t) =f(x_0,t),~
\alpha \in (0,1],
\end{array}
\right.
\end{equation}
with appropriate initial conditions on $x_i(0)$ given. Here
$p=\lceil \alpha \rceil$. For simplicity we assume $x_i(t) \in
H^1(\Omega,\mathcal{T})$ except $x_m(t) \in
L^2(\Omega,\mathcal{T})$. We will continue to use the upwind fluxes
to seek $X_i$, such that for all $v_i  \in V$, the following holds

\begin{equation}\label{Scheme4}
\left\{\begin{array}{l}
\displaystyle\big(X_{i+1},v_i\big)_{I_j}+\bigg(X_i,\frac{dv_i}{dt}\bigg)_{I_j}-\big(
X_i(t_j^-)v_i(t_j^-)-X_i(t_{j-1}^-)v_i(t_{j-1}^+)\big)=0, \\
\qquad i=0, \ldots, \max(m,p)-1, X_0(t) = x(t); \\
\\
\displaystyle \big( {_0D_t^{-(p-\alpha)}}X_p,v_m
\big)_{I_j}+\big(X_m,v_m\big)_{I_j}=\big(f(X_0,t),v_m \big)_{I_j},
\end{array}
\right.
\end{equation}
subject to the appropriate initial conditions.

The analysis of this scheme is generally similar to that of the
previous one and will not be discussed further, although, as we
shall illustrate shortly, there are details that remain open. A main
difference is that the order of super convergence at the downwind
point changes to $k+1+\min\{k,\max\{\alpha,m\}\}$ and the impact of
the fractional operator is thus eliminated by the linear classic
operator as long as $m\geq \alpha$. However, for the case where $
\lfloor \alpha \rfloor \geq k,m$, the situation is less clear.

Let us first consider  a linear example to illustrate that the order
of super convergence  $2k+1$ at downwind points can also be obtained
when $\alpha$ is an integer at the left end point of the interval,
provided the initial condition is not overspecified. On the
computational domain $t\in \Omega=(0,1)$, we consider
\begin{equation} \label{ExampleL1Prime}
{_0^CD_t^\alpha
x(t)}+\frac{dx(t)}{dt}=-2x(t)+\frac{\Gamma(6)}{\Gamma(6-\alpha)}t^{5-\alpha}+2t^5+5t^4+2,~~~
\alpha \in [0,1],
\end{equation}
with the initial condition $x(0)=1$ and the exact solution
$x(t)=t^5+1$.

\begin{figure}
\vbox{ \hbox{
\includegraphics[width=2.5in,angle=0]{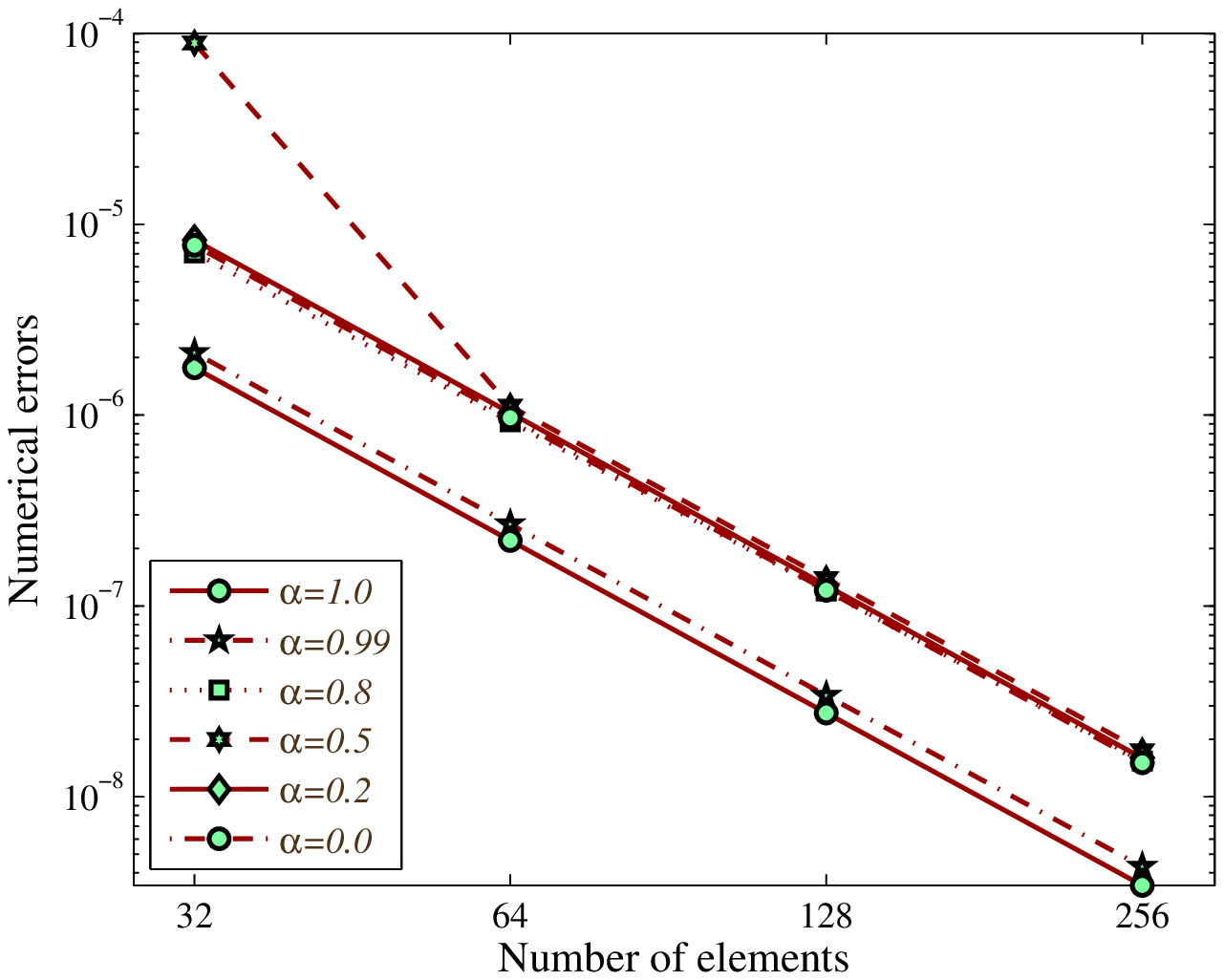}
\includegraphics[width=2.5in,angle=0]{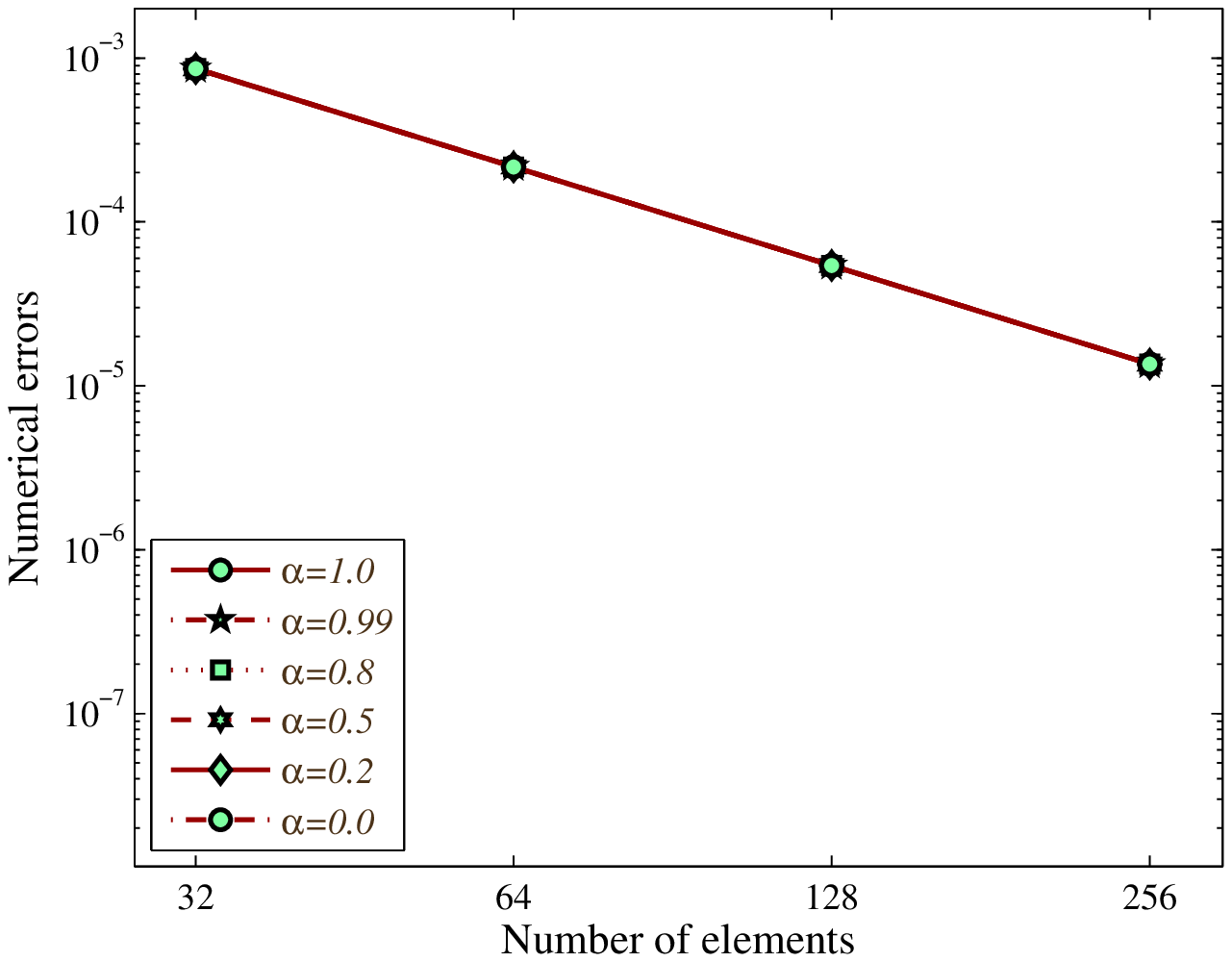}}
\hbox{
\includegraphics[width=2.5in,angle=0]{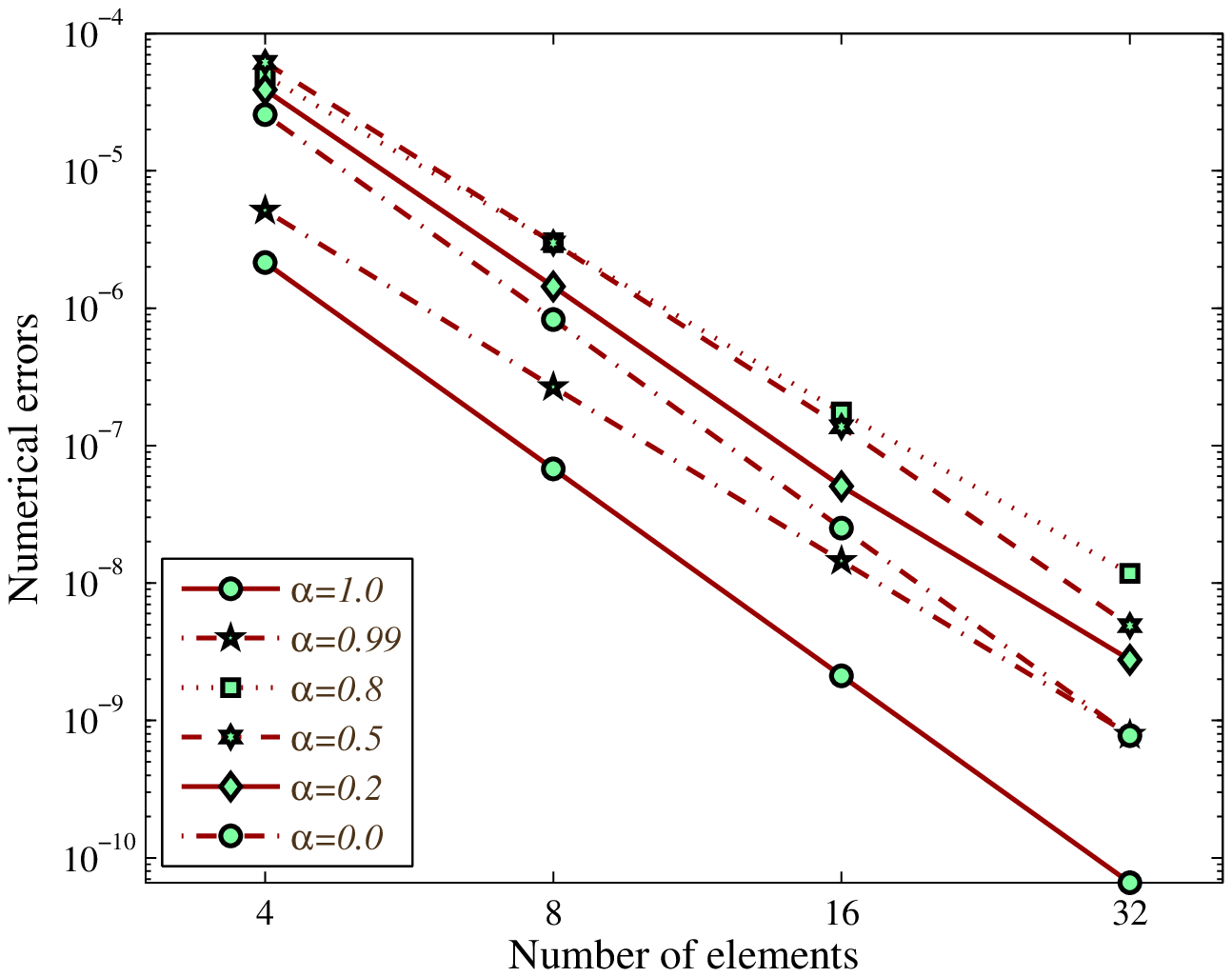}
\includegraphics[width=2.5in,angle=0]{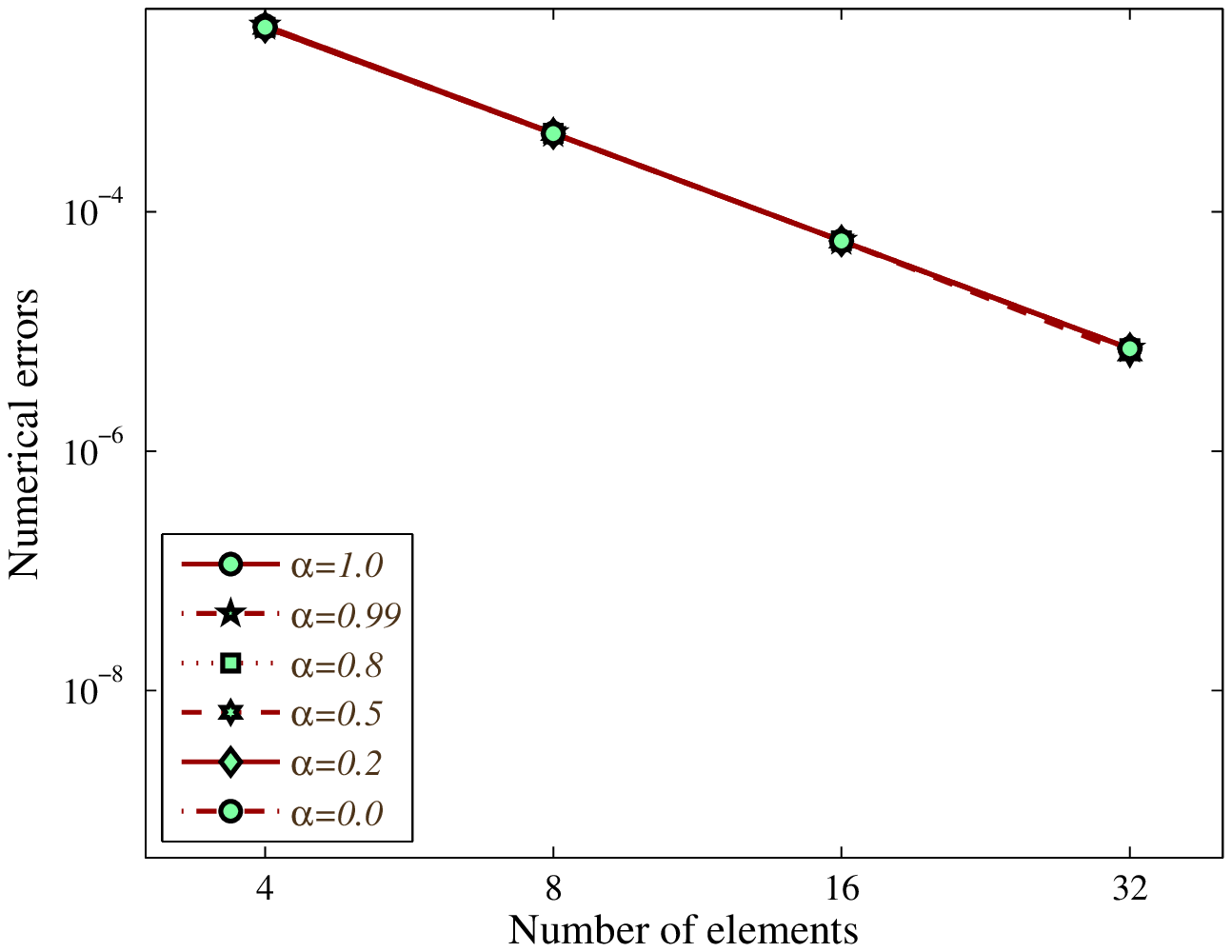}}
\hbox{
\includegraphics[width=2.5in,angle=0]{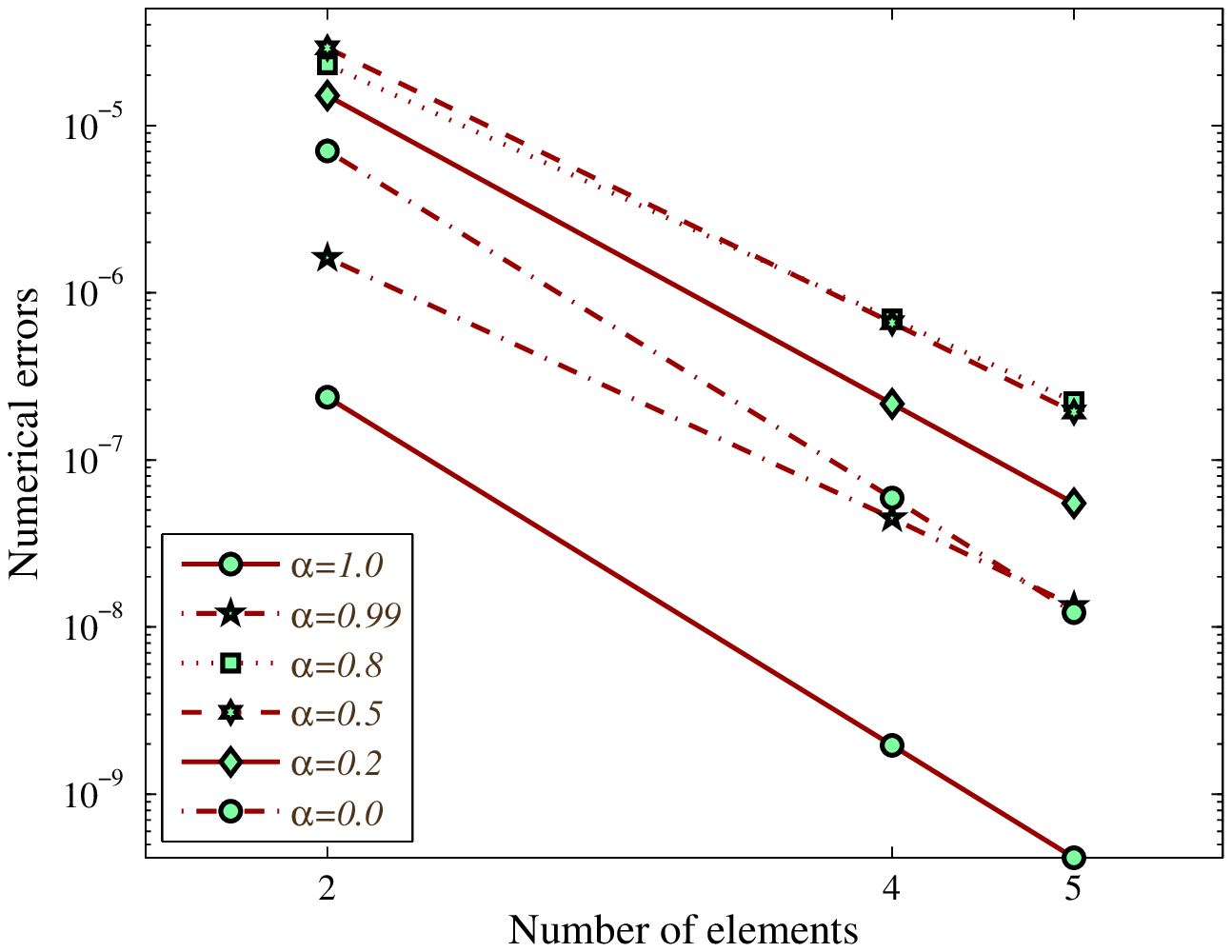}
\includegraphics[width=2.5in,angle=0]{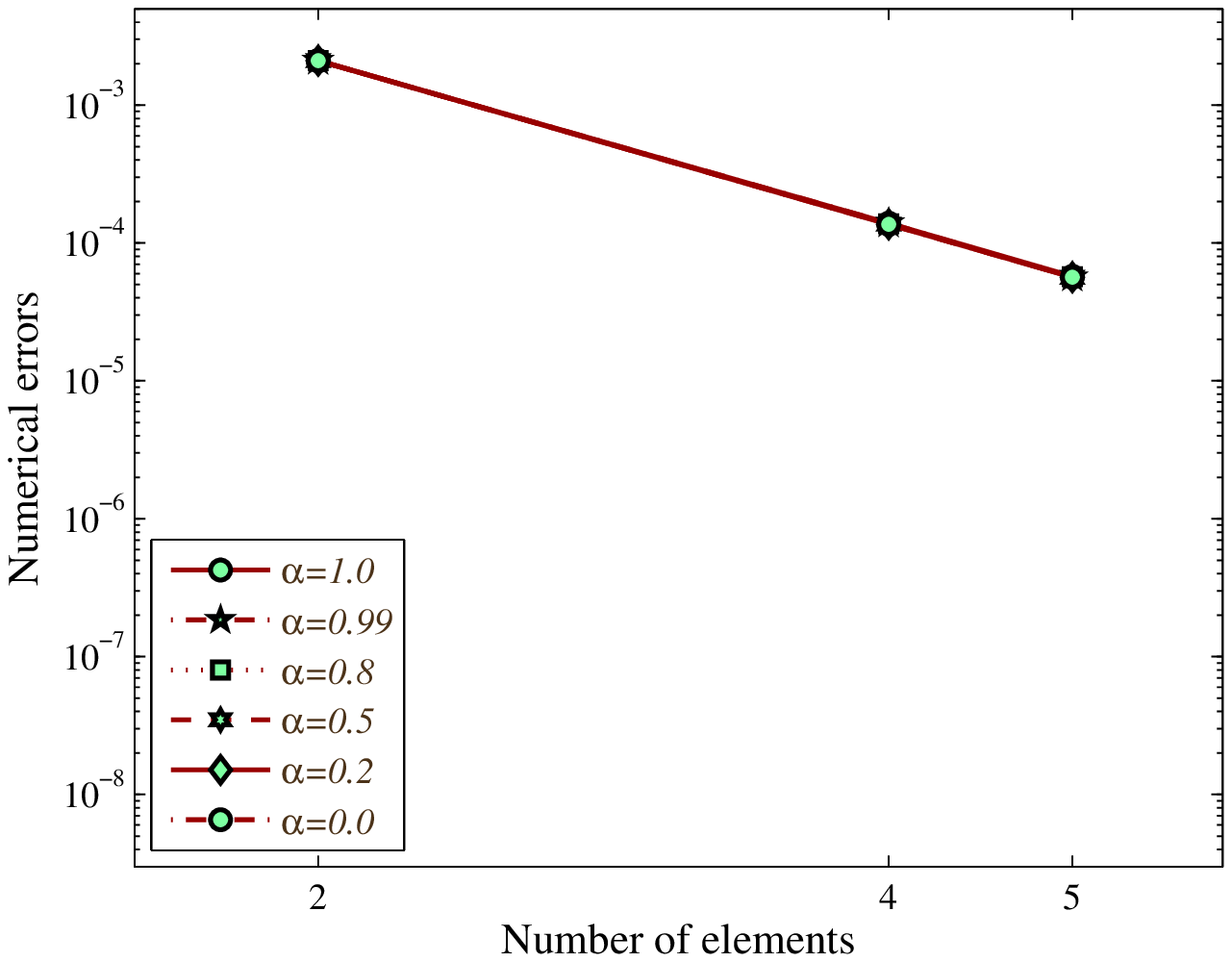}}}
\caption{The convergence of (\ref{ExampleL1Prime}) for $k=1$ (top
row), $k=2$ (middle row) and $k=3$ (bottom row). In the left column
we show the convergence at the downwind points while the right
column displays $L^2$-convergence. \label{L1Prime}}
\end{figure}
As expected, Fig. \ref{L1Prime} confirms super convergence of $k+1+\min\{k,\max\{\alpha,m\}\}$ at the downstream point.

\begin{figure}
\vbox{ \hbox{
\includegraphics[width=2.5in,angle=0]{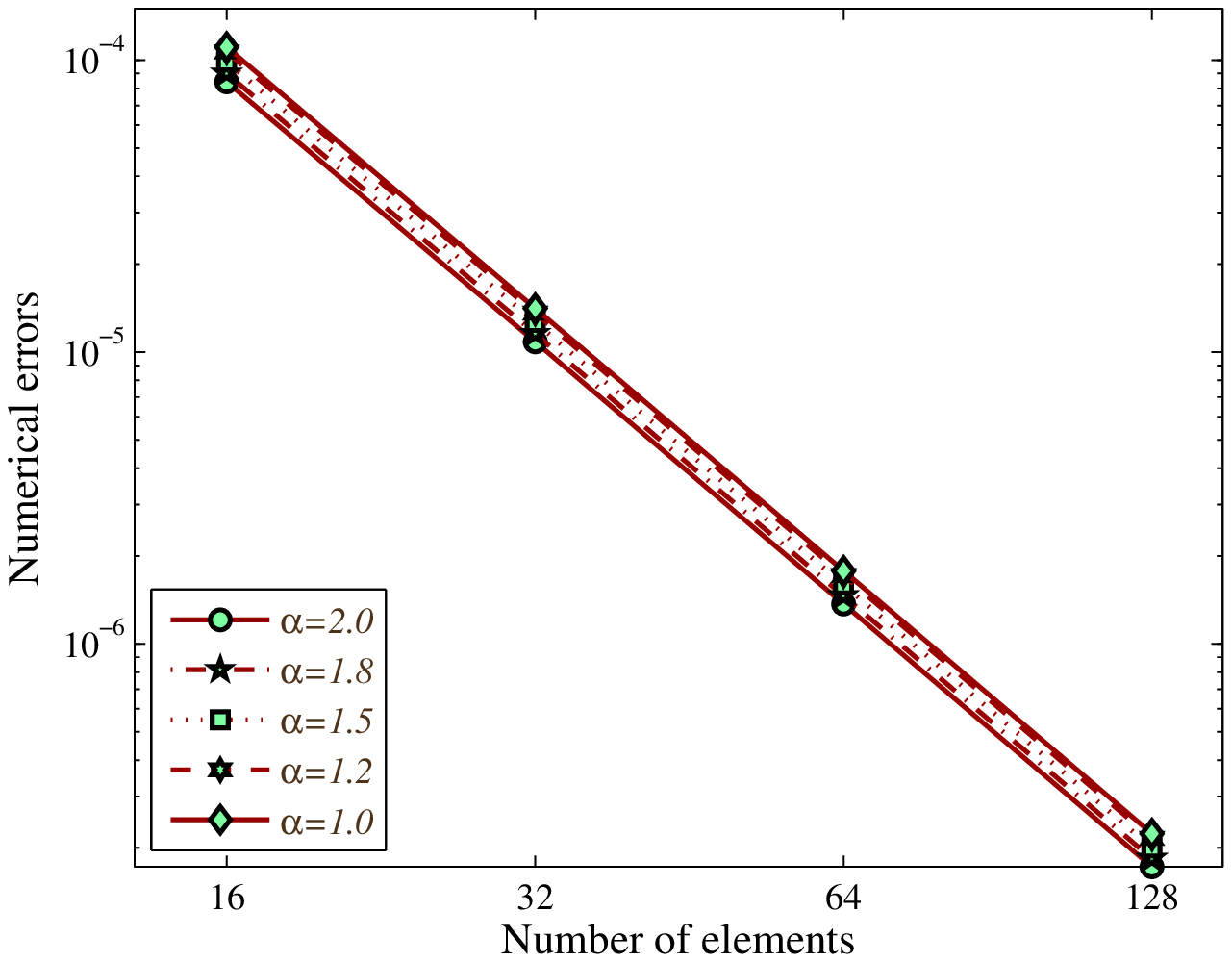}
\includegraphics[width=2.5in,angle=0]{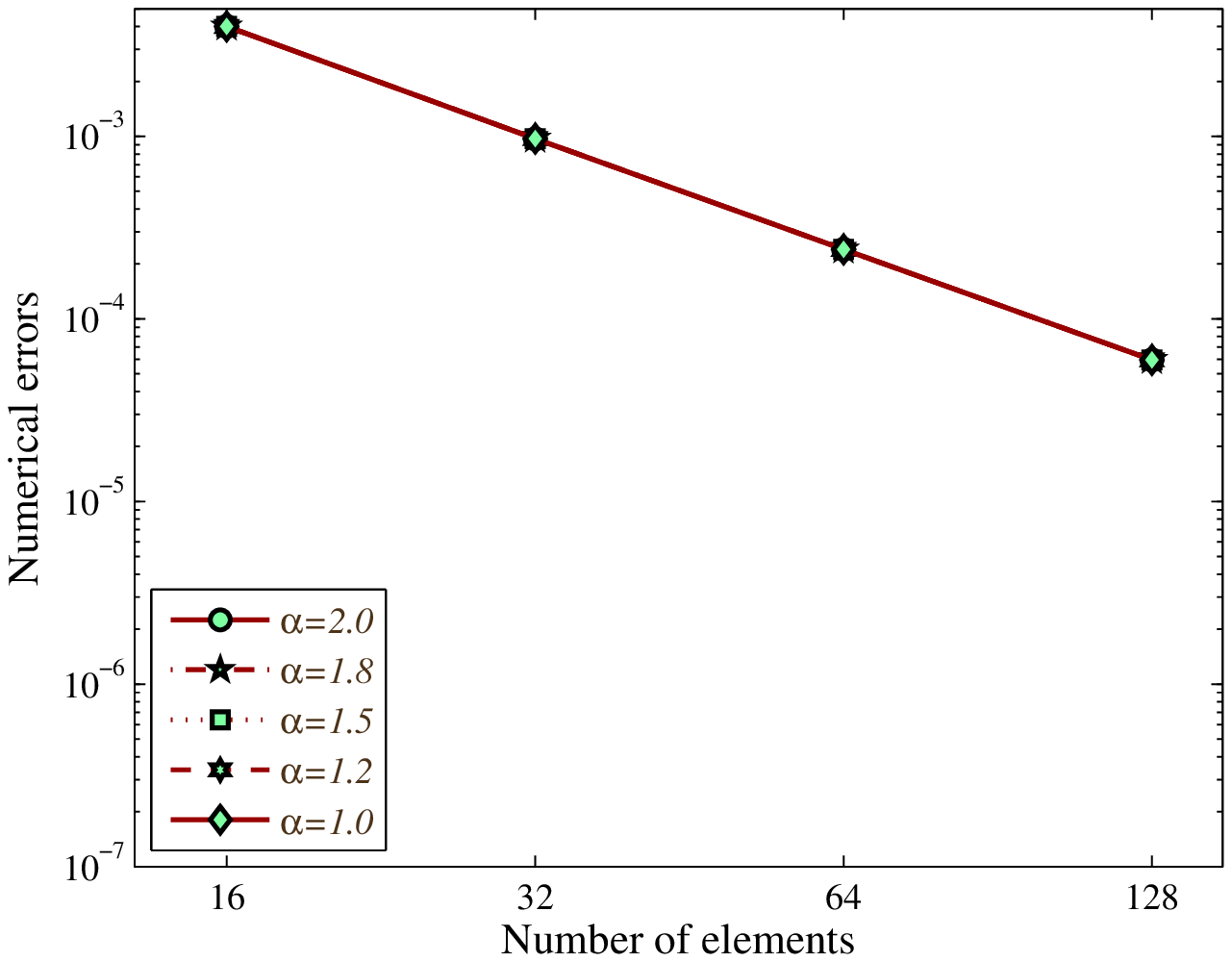}}
\hbox{
\includegraphics[width=2.5in,angle=0]{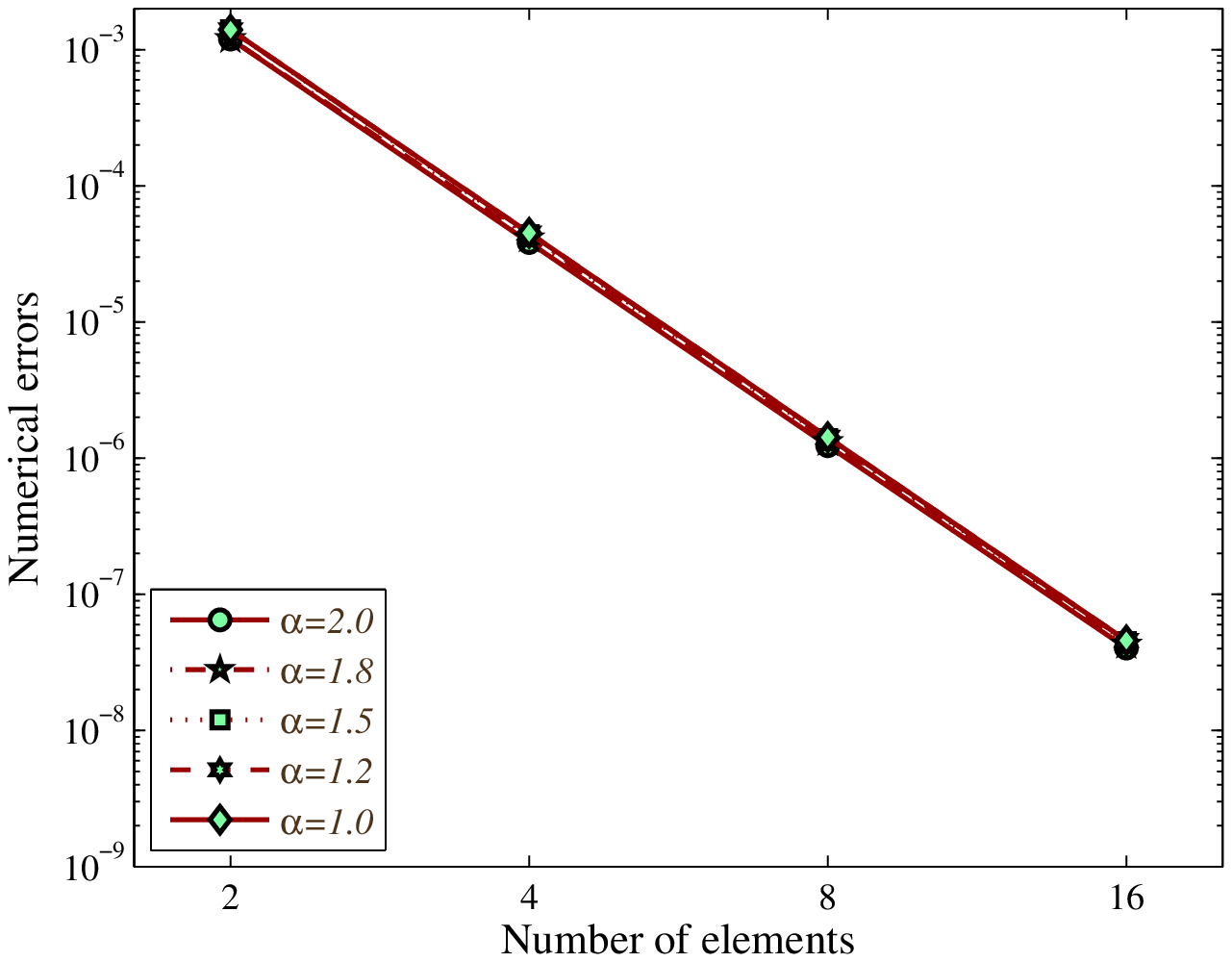}
\includegraphics[width=2.5in,angle=0]{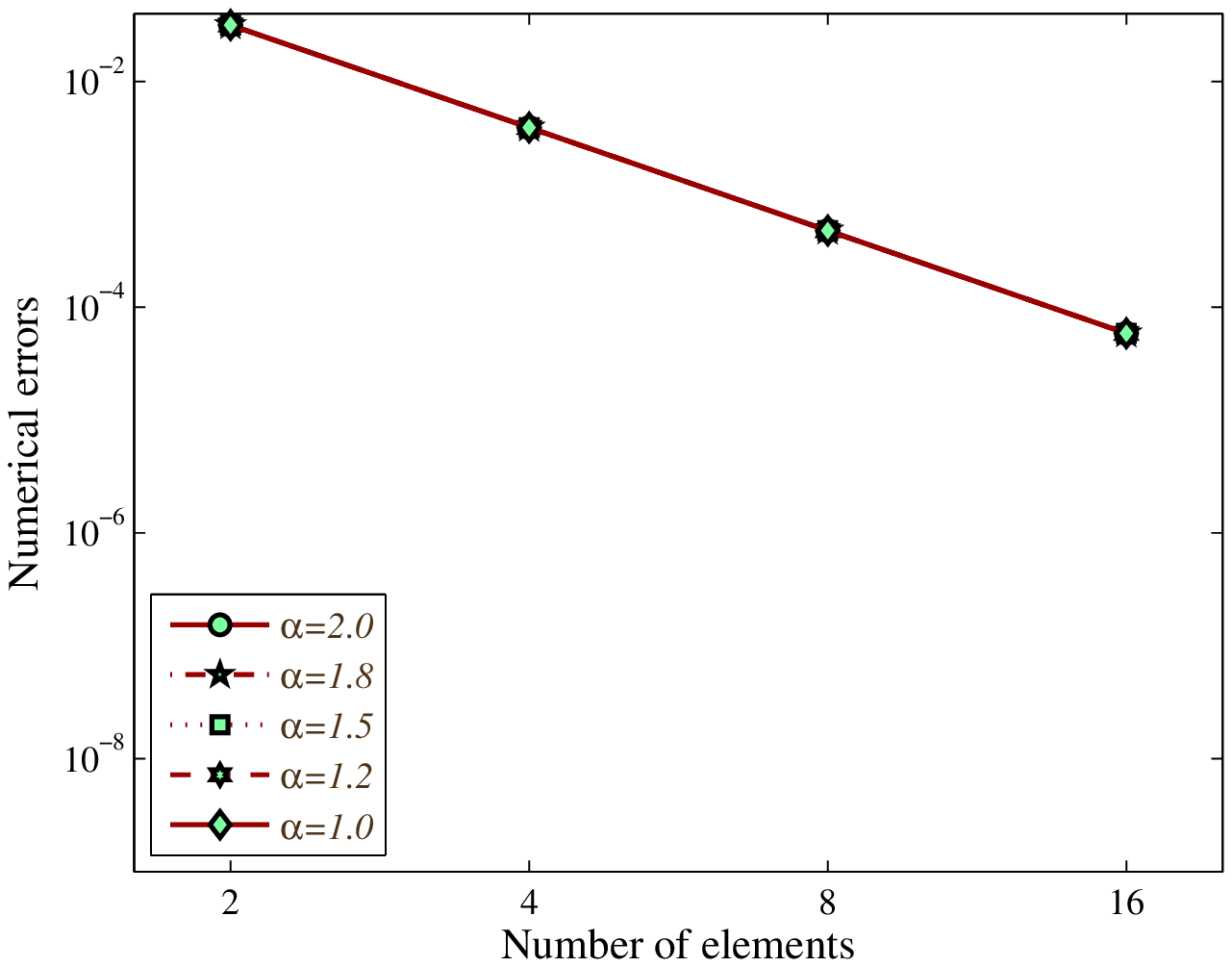}}
\hbox{
\includegraphics[width=2.5in,angle=0]{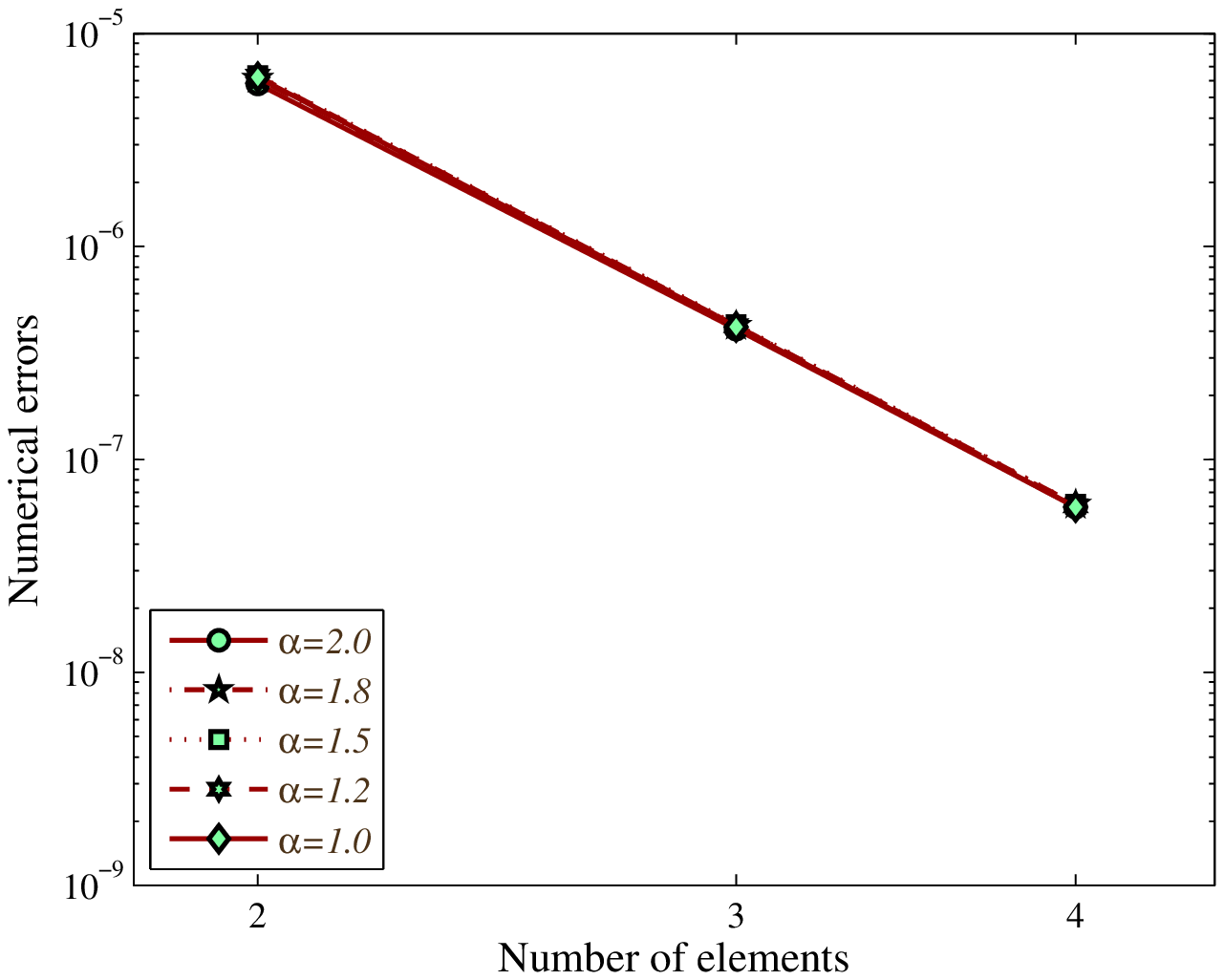}
\includegraphics[width=2.5in,angle=0]{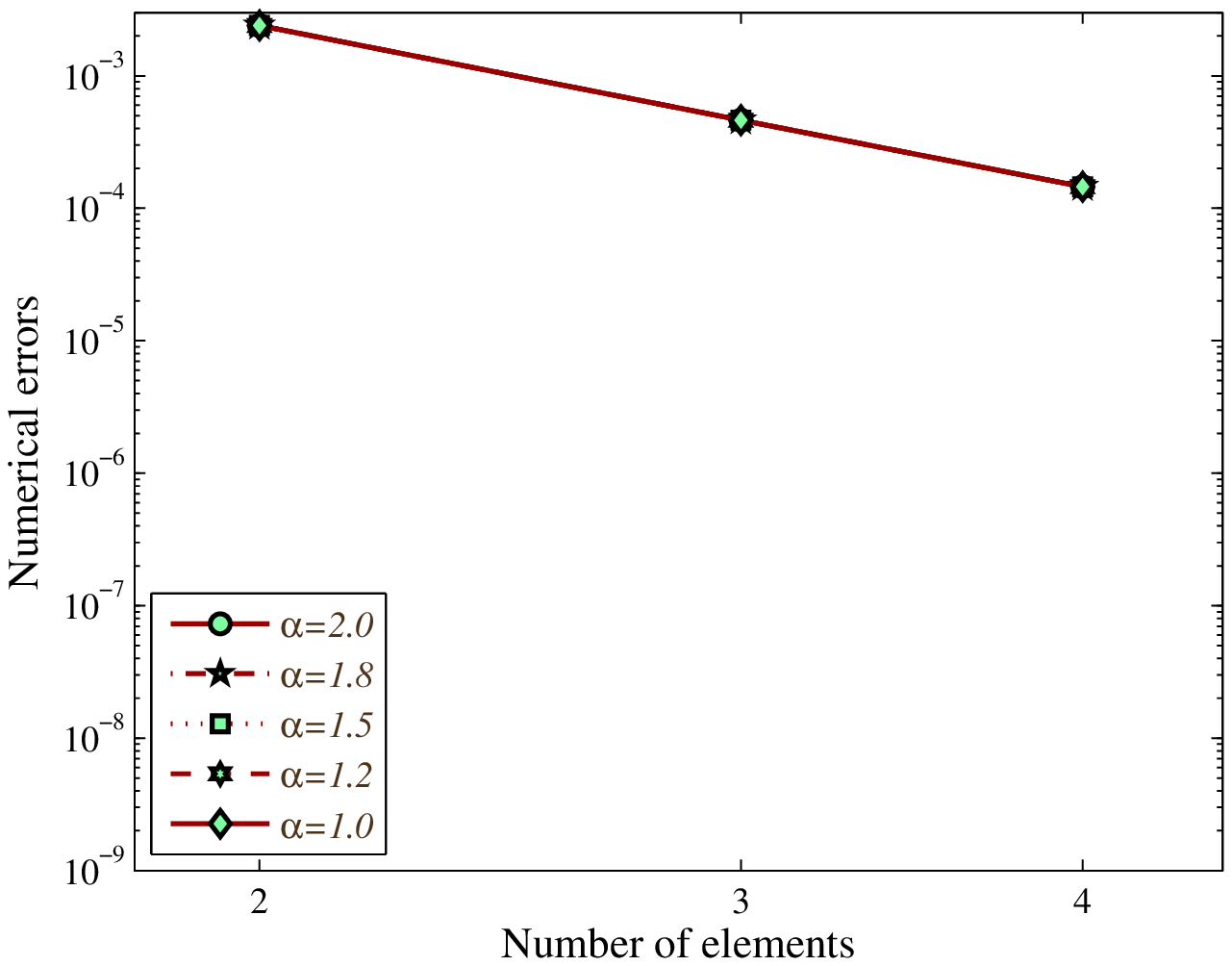}}}
\caption{The convergence of (\ref{Scheme4}) for $k=1$ (top row),
$k=2$ (middle row) and $k=3$ (bottom row) in (\ref{ExampleN5}). In
the left column we show the convergence at the downwind points while
the right column displays $L^2$-convergence. \label{N5}}
\end{figure}

Let us  consider a nonlinear problem, given for $t\in \Omega=(0,1)$,
as
\begin{equation} \label{ExampleN5}
~~~~\begin{array}{lll} \displaystyle{_0^CD_t^\alpha
x(t)}+\frac{d^3x(t)}{dt^3} &=& \displaystyle
-2x^2(t)+\frac{\Gamma(6)}{\Gamma(6-\alpha)}t^{5-\alpha}+\frac{\Gamma(1)}{\Gamma(3-\alpha)}t^{2-\alpha}+2t^{10}+2t^7
\\
&& \displaystyle+4t^6 +4t^5+0.5t^4+2t^3+64t^2+4t+2,~~ \alpha \in
[1,2],
\end{array}
\end{equation}
with the initial condition $x(0)=1$, $x^\prime(0)=1$,
$x^{\prime\prime}(0)=1$ and the exact solution
$x(t)=t^5+\frac{1}{2}t^2+t+1$. We note that in this case, $\alpha
\in [1,2]$ but $m=3$ and, as shown in Fig. \ref{N5} super
convergence of order $2k+1$ is maintained in this case.

\begin{figure}
\vbox{ \hbox{
\includegraphics[width=2.5in,angle=0]{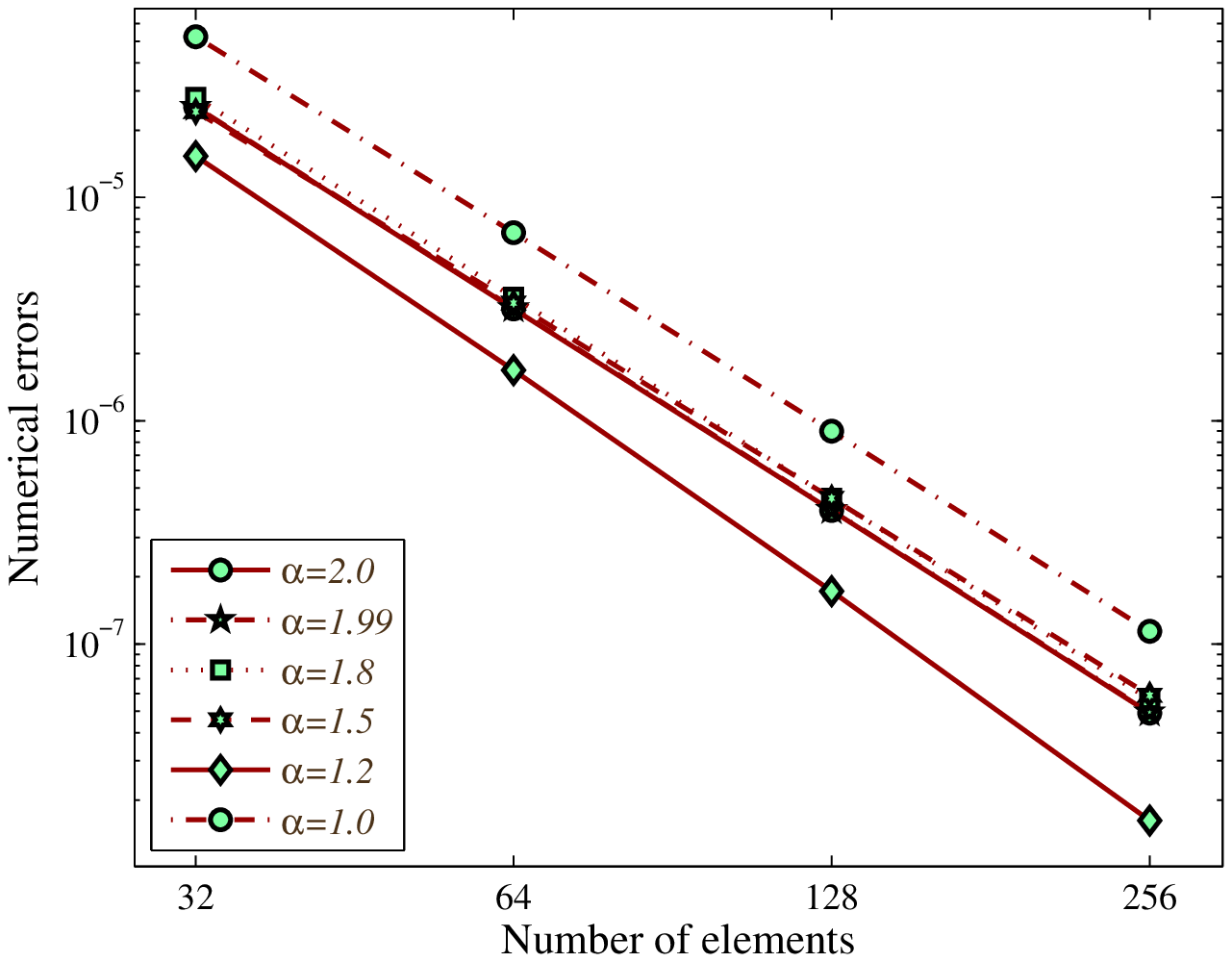}
\includegraphics[width=2.5in,angle=0]{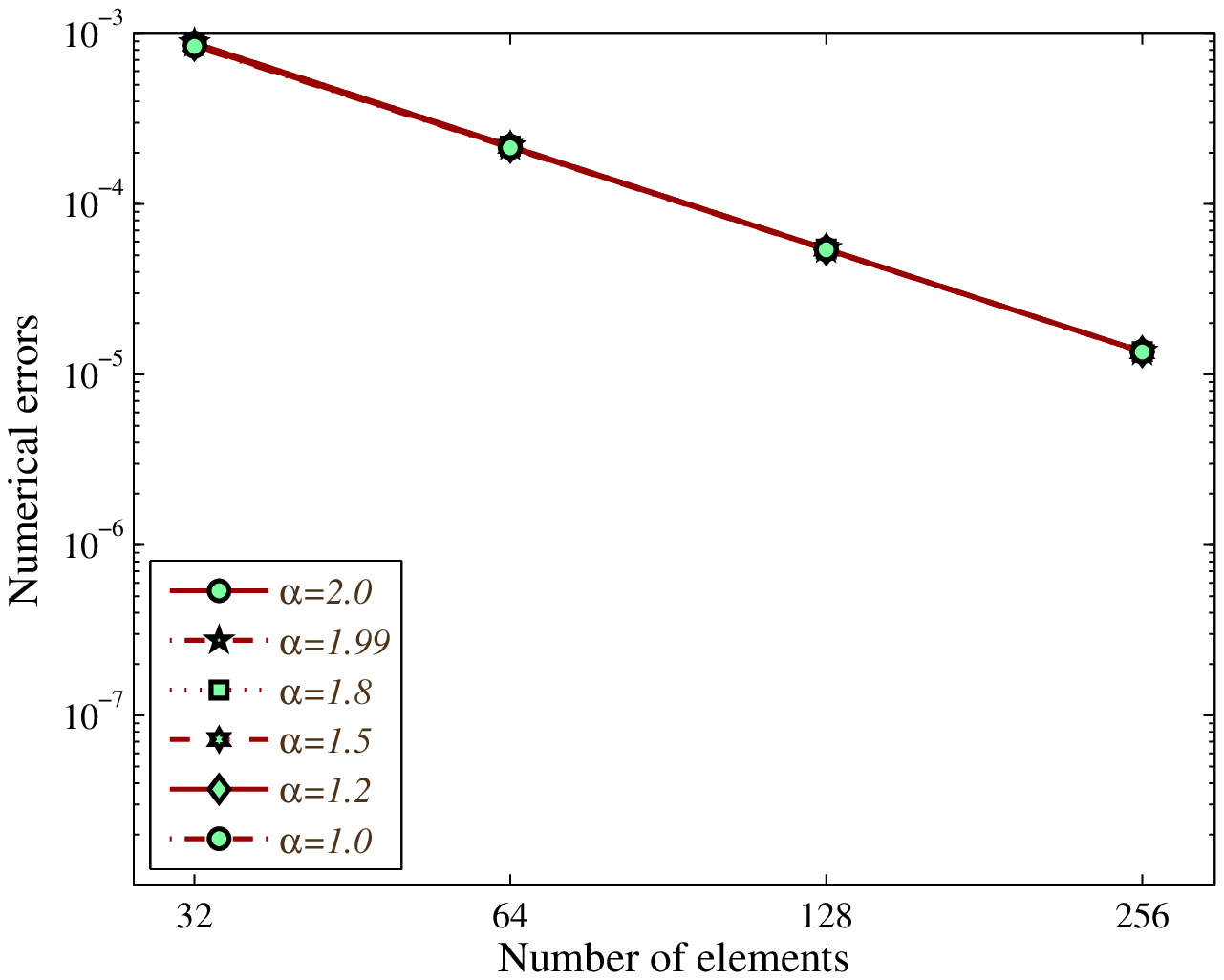}}
\hbox{
\includegraphics[width=2.5in,angle=0]{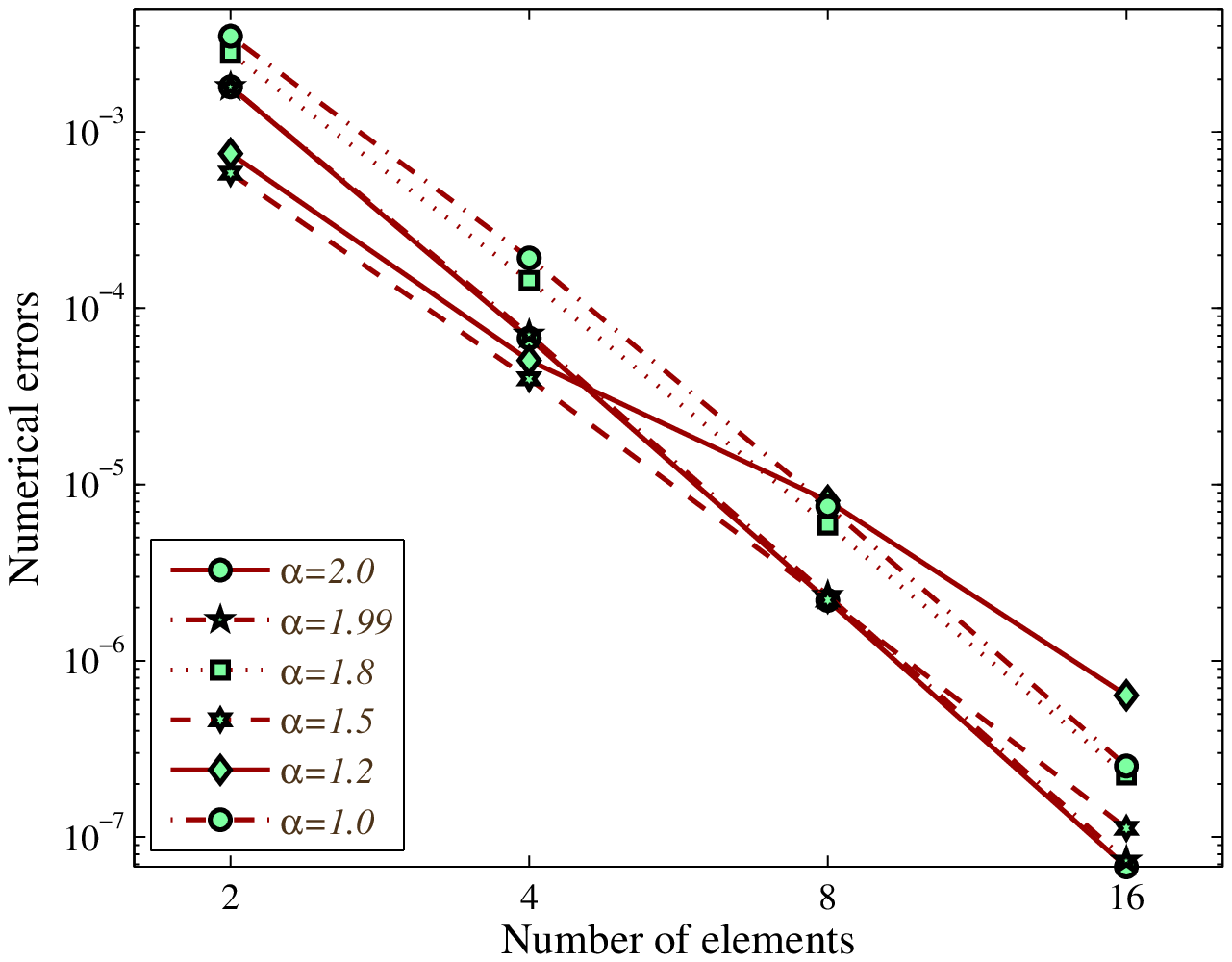}
\includegraphics[width=2.5in,angle=0]{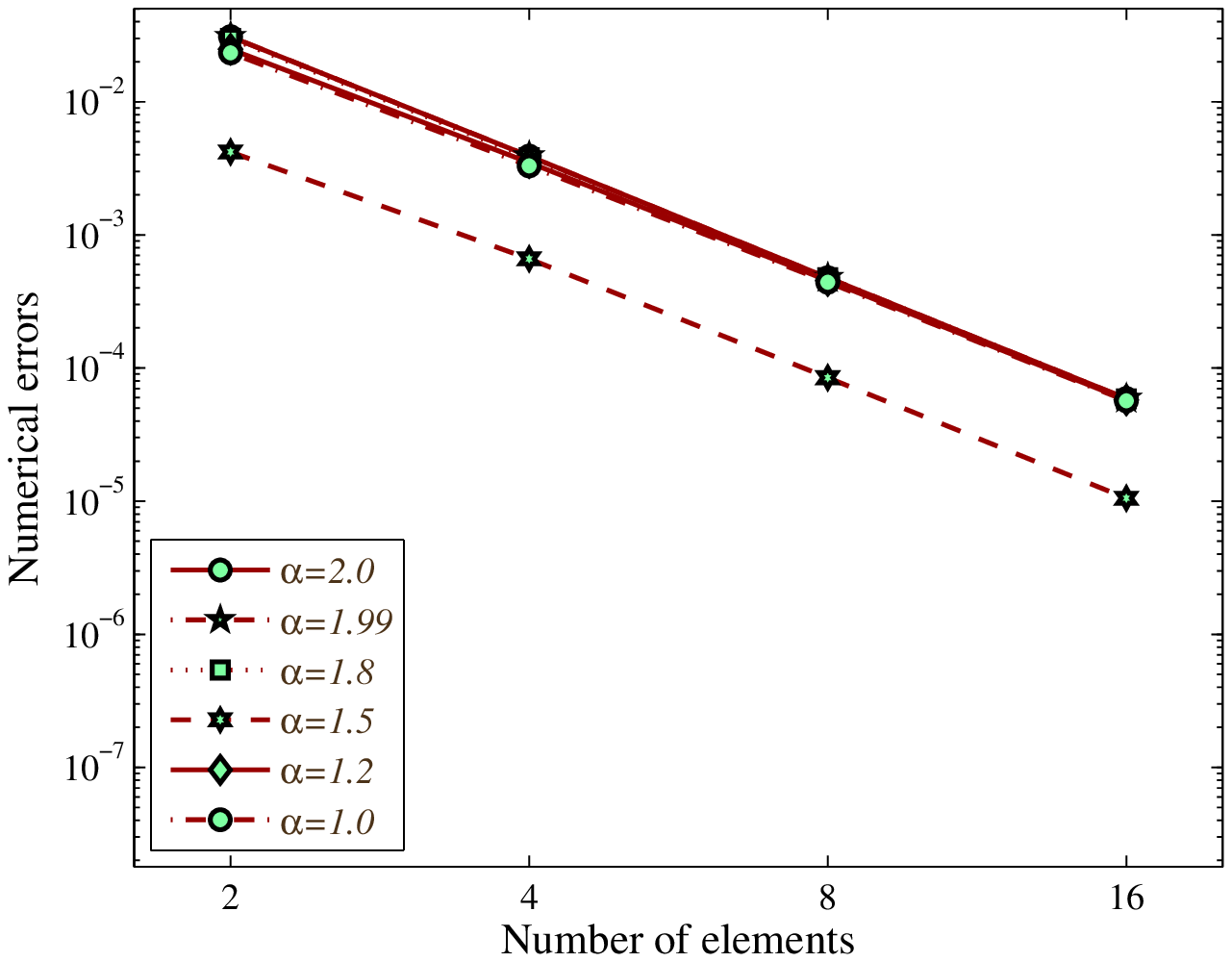}}
\hbox{
\includegraphics[width=2.5in,angle=0]{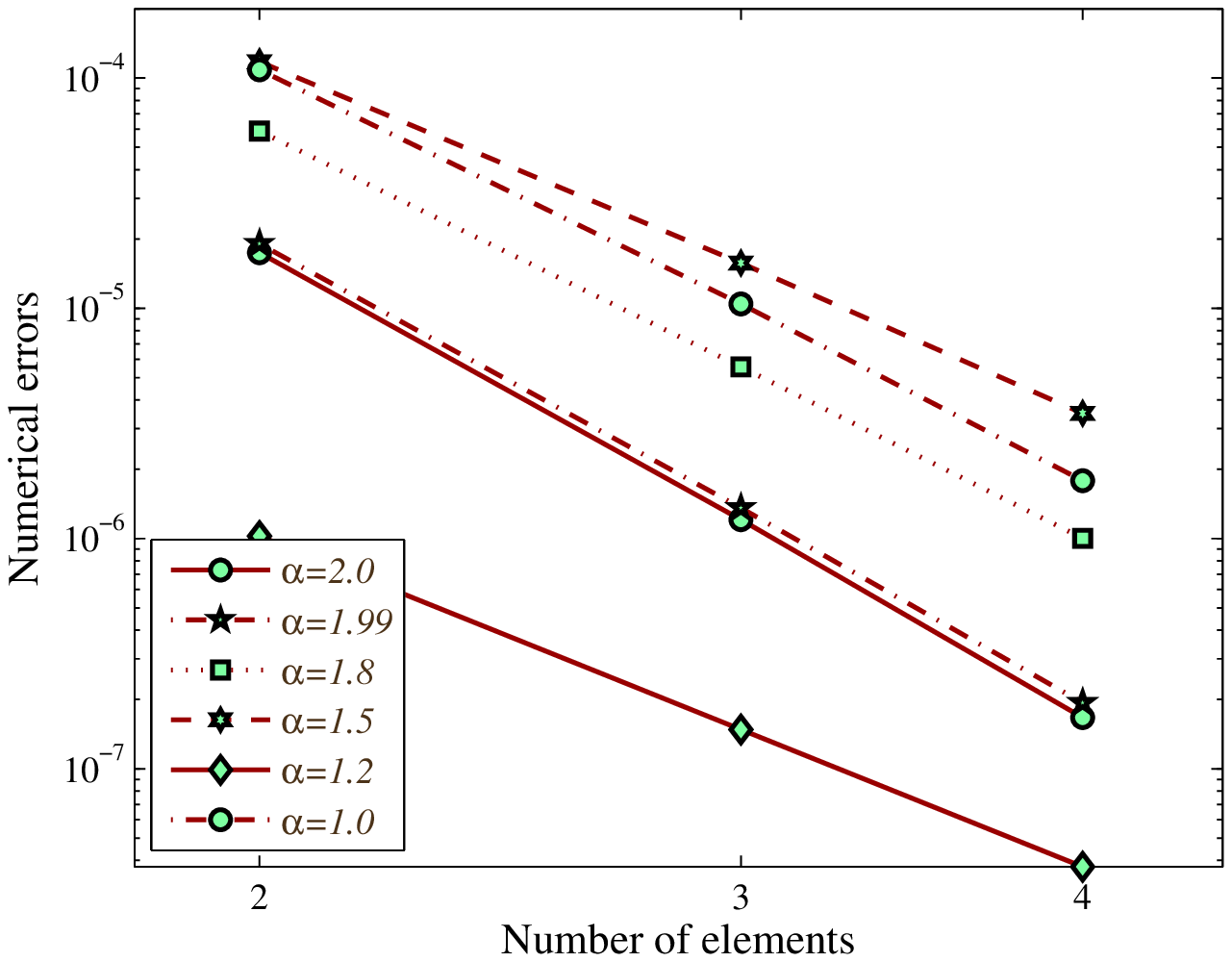}
\includegraphics[width=2.5in,angle=0]{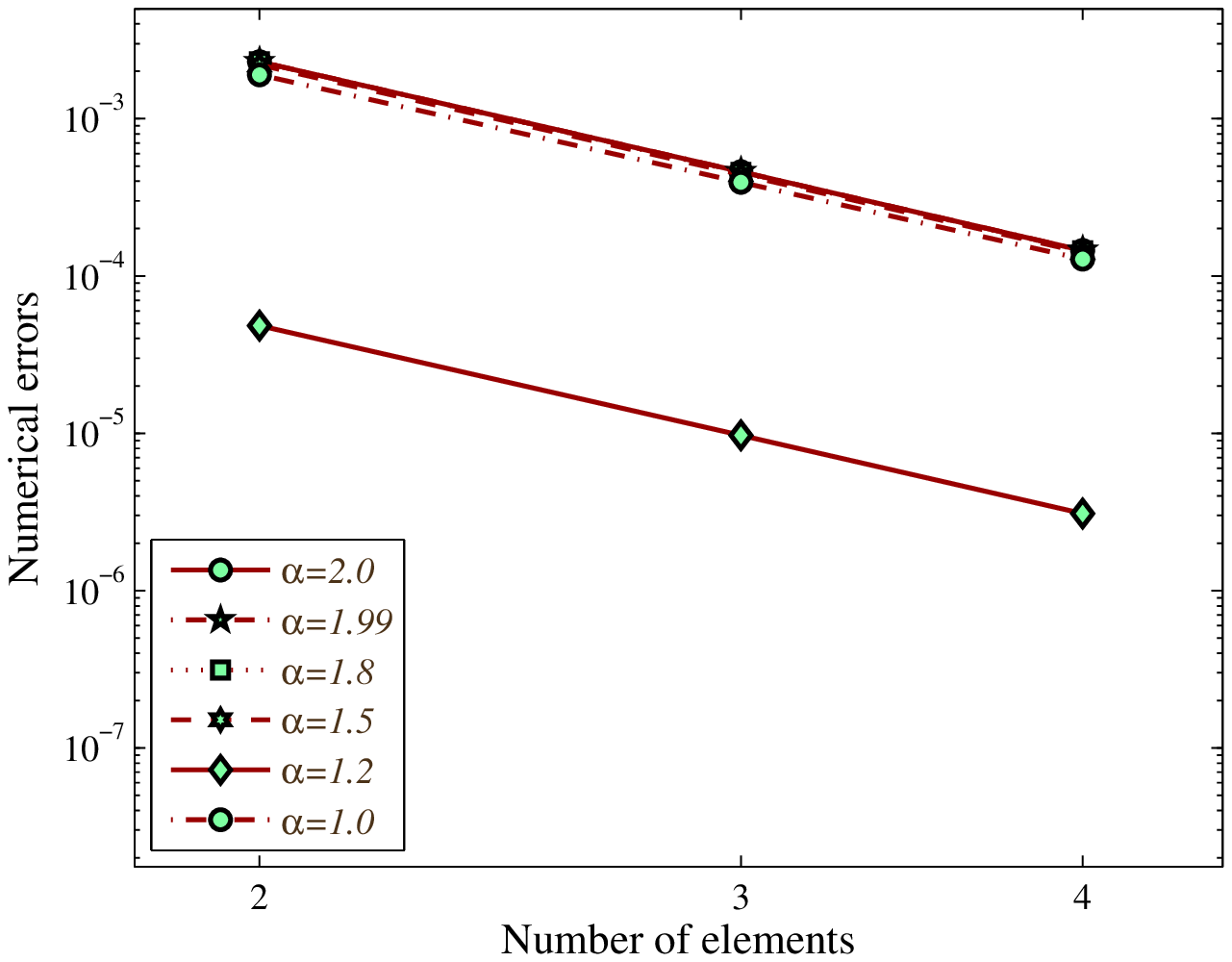}}}
\caption{The convergence of (\ref{Scheme4}) for $k=1$ (top row),
$k=2$ (middle row) and $k=3$ (bottom row) in (\ref{ExampleN4}). In
the left column we show the convergence at the downwind points while
the right column displays $L^2$-convergence. \label{N4}}
\end{figure}

However, if the assumption that $m\geq \lceil \alpha \rceil $ is
violated, an $\alpha$-dependent rate of convergence re-emerges as
illustrated in the following example. Consider a nonlinear problem,
given for $t\in[0,1]$ as

\begin{equation}
\label{ExampleN4} ~~~~~~{_0^CD_t^\alpha
x(t)}+\frac{dx(t)}{dt}=-2x^2(t)+\frac{\Gamma(6)}{\Gamma(6-\alpha)}t^{5-\alpha}+5t^4+1+2(t^5+t+1)^2,~
\alpha \in [1,2],
\end{equation}
with the initial condition $x(0)=1$, $x^\prime(0)=1$ and the exact
solution $x(t)=t^5+t+1$. Figure \ref{N4} shows superconvergence of
$k+1+\min\{k,\alpha\}$ at the downstream point.

As a final example, let us consider \begin{equation}
\label{ExampleL2} {_0^CD_t^\alpha
x(t)}=-2x(t)+\frac{\Gamma(6)}{\Gamma(6-\alpha)}t^{5-\alpha}+2t^5+2t+2,~~~
\alpha \in [1,2],
\end{equation}
on the computational domain $t\in \Omega=(0,1)$, with the initial
condition $x(0)=1$, $x^\prime(0)=1$ and the exact solution
$x(t)=t^5+t+1$.

\begin{figure}
\hbox{
\includegraphics[width=2.5in,angle=0]{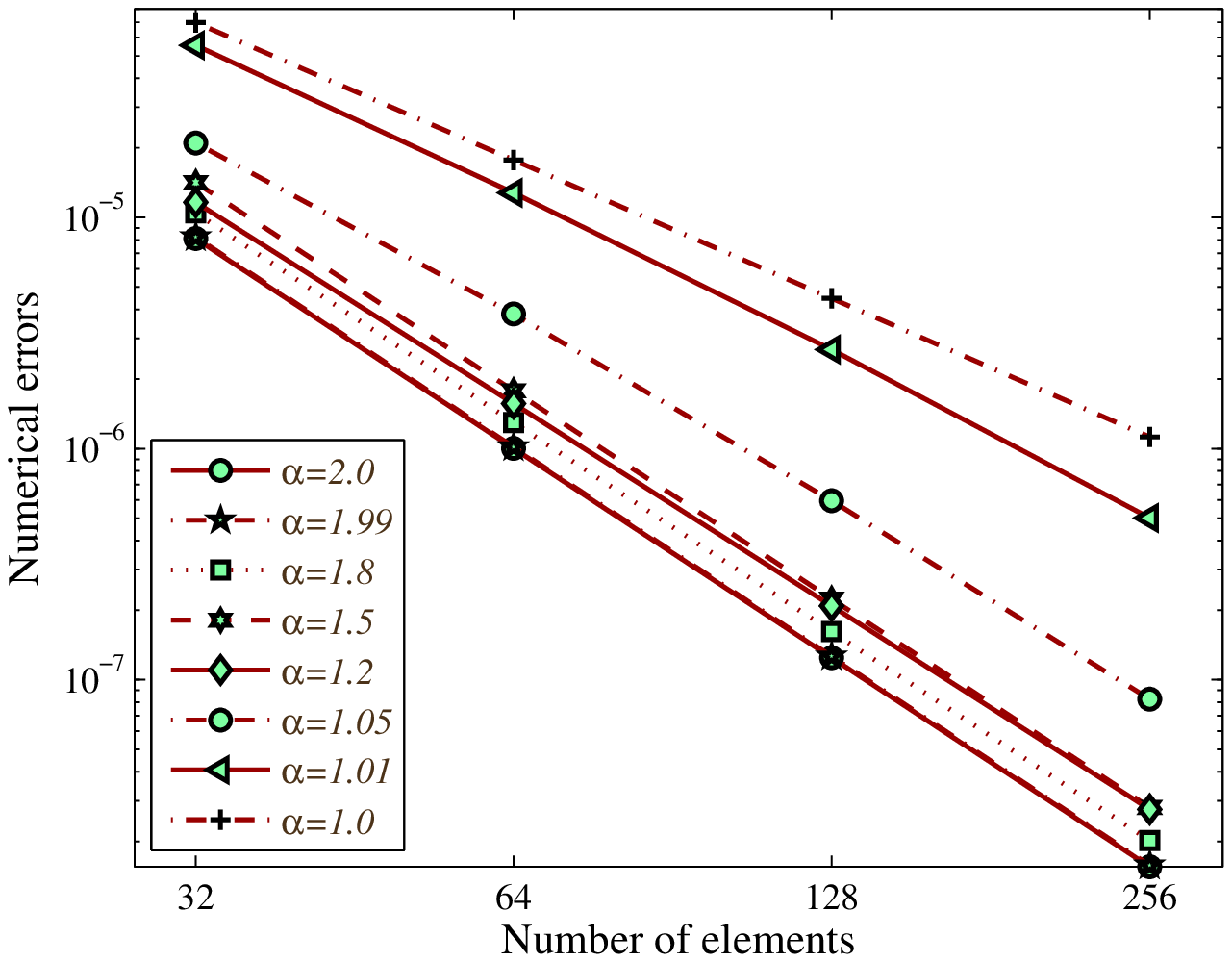}
\includegraphics[width=2.5in,angle=0]{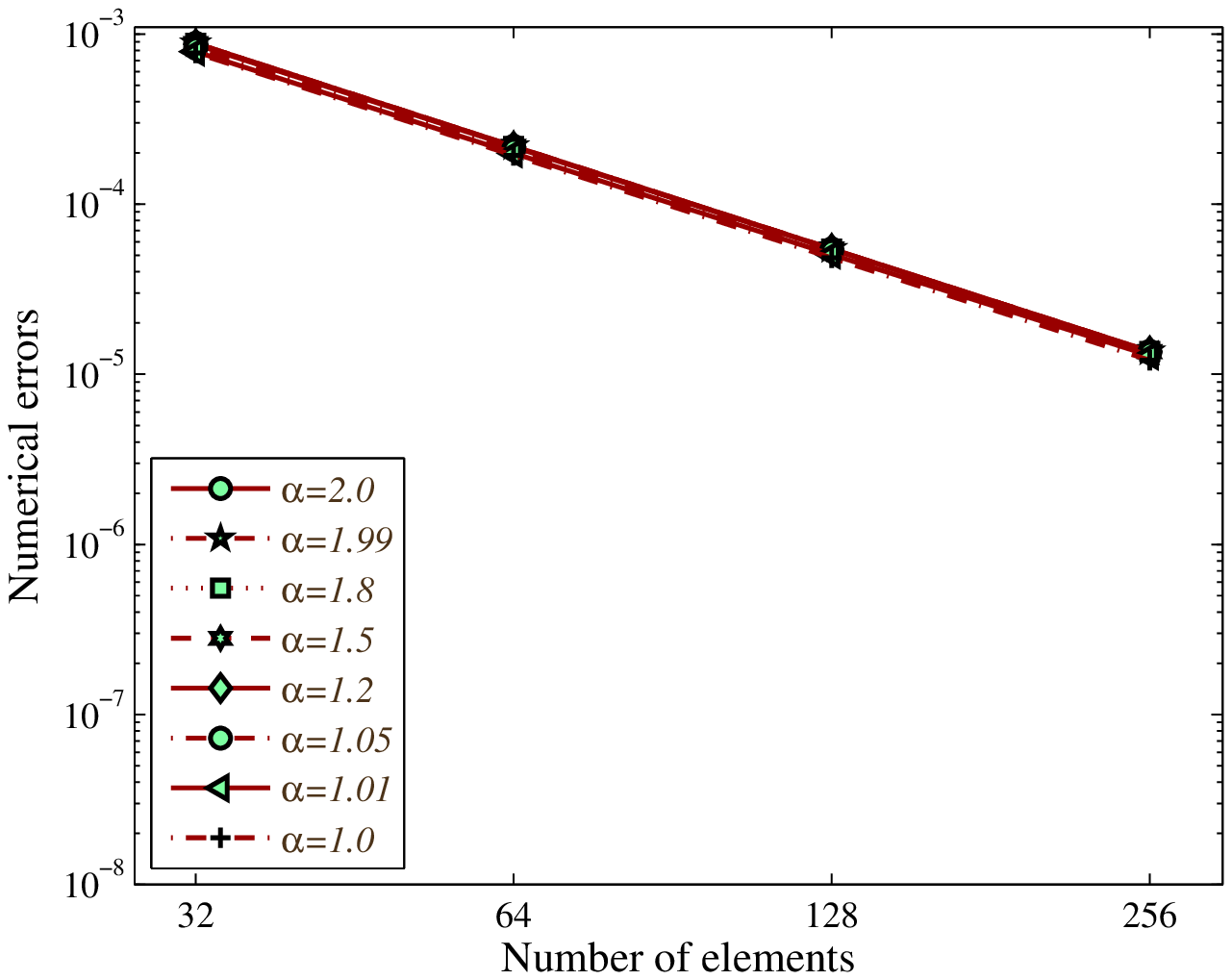}}
\caption{The convergence of (\ref{Scheme4}) for $k=1$ in
(\ref{ExampleL2}). On the left we show the convergence at the
downwind points while the right figure displays $L^2$-convergence.
\label{L2a}}
\end{figure}

In Fig. \ref{L2a} we show the results for $k=1$. Following the
previous analysis, we would expect an order of convergence as
$k+1+\min\{k,\alpha\}$ which in this case would be third order.
However, the results in Fig. \ref{L2a} highlights a reduction in the
order of convergence at the endpoint as $\alpha$ approaches the
value one. The mechanism for this is not fully understood but is
likely associated with an over specification of the initial
conditions in this singular limit. Increasing $k$ recovers the
expected convergence rate for all values of $\alpha$.

\section{Concluding remarks}
We introduce an LDG schemes with upwind fluxes for general FODEs.
The schemes enable an element by element solution, hence avoiding
the need for a full global solve. Through analysis, we highlight
that the  scheme converges with the optimal order of convergence
order $k+1$ in $L^2$ norm and shows  superconvergence  at the
downwind point of each interval with an order of   convergence order
of $k+1+\min\{k,\alpha\}$, where $\alpha$ refers to the order of the
fractional derivatives and $k$ the degree of the approximating
polynomial. We discuss the mechanism for this superconvergence and
extend the discussion to cases where classic ODE terms of order $m$
are included in the equation. In this case, the order of the super
convergence becomes $k+1+\min\{k,\max\{\alpha,m\}\}$, i.e., the
behavior of the classic derivative dominates that of the fractional
operator provided $m\geq  \lceil \alpha \rceil $. This is confirmed
through examples. A final case in which $k\leq \lfloor \alpha
\rfloor$ shows that in this case, the expected convergence of
$k+1+\min\{k,\alpha\}$ is violated as $\alpha$ approaches one. The
mechanism for this remains unknown and we hope to report on that in
future work.


\begin{thebibliography}{}





\bibitem{Adjerid:02}
Adjerid, S.,  Devine, K.D., Flaherty, J.E.,  Krivodonova, L.: A
posteriori error estimation for discontinuous Galerkin solutions of
hyperbolic problems. Comput. Methods Appl. Mech. Engrg. {\bf191},
 1097--1112 (2002)

\bibitem{Bassi:97} Bassi, F.,   Rebay, S.: A high-order accurate discontinuous finite element method for the numerical solution
of the compressible Navier-Stokes equations. J. Comput. Phys.
{\bf131}, 267-279 (1997)

\bibitem{Brunner:06}
 Brunner, H.,  Sch\"tzau, D.: Hp-Discontinuous Galerkin time-stepping
for Volterra integrodifferential equations.
 SIAM J. Numer. Anal. {\bf44}, 224-245 (2006)

\bibitem{Butzer:00} Butzer, P.L.,   Westphal,  U.: An Introduction to Fractional Calculus. World
Scientific, Singapore (2000)

\bibitem{Cockburn:03}  Cockburn,  B.: Discontinuous Galerkin method. Z. Angew. Math. Mech.
{\bf83},
731-754 (2003)

\bibitem{Cockburn:98}  Cockburn, B.,   Shu, C.-W.: The local discontinuous Galerkin method for time-dependent
convection diffusion systems. SIAM J. Numer. Anal. {\bf35},
2440-2463 (1998)

\bibitem{Delfour:81}  Delfour, M.,  Hager, W.,  Trochu,  F.:
 Discontinuous Galerkin methods for ordinary differential equations.  Math. Comp. {\bf36},
 455-473 (1981)

\bibitem{Deng:07}  Deng,  W.H.: Numerical algorithm for the time fractional Fokker-Planck equation. J. Comput. Phys. {\bf227},  1510-1522 (2007)

\bibitem{Deng:11} Deng  W.H.,  Hesthaven, J.S.: Discontinuous Galerkin methods for fractional diffusion
equations. ESAIM: M2AN {\bf47},  1845-1864 (2013)

\bibitem{Diethelm:02}  Diethelm, K.,  Ford, N.J.,  Freed, A.D.: A predictor-corrector approach
for the numerical solution of fractional differential equations.
Nonl. Dynam. {\bf29},  3-22 (2002)

\bibitem{Guo:06} Guo B.Q.,  Heuer, N.: The optimal convergence of the h-p version of the boundary element method with
quasiuniform meshes for elliptic problems on polygonal domains. Adv.
Comput. Math. {\bf24},  353-374 (2006)

\bibitem{Hartley:08}  Hartley, T.T., Lorenzo, C.F.,   Qammer, H.K.:
 Chaos in a fractional order Chua's system. IEEE Trans. Circuits Syst. I: Fundam. Theory Appl. {\bf42},
 485-490 (1995)

\bibitem{Hesthaven:04}  Hesthaven, J.S.,  Warburton, T.:
 High-order nodal discontinuous Galerkin methods for Maxwell eigenvalue problem. Royal Soc. London Ser. A {\bf362},
 493-524 (2004)

\bibitem{Hesthaven:08}  Hesthaven, J.S.,   Warburton, T.:  Nodal Discontinuous Galerkin Methods: Algorithms, Analysis, and Applications.
 Springer-Verlag, New York, USA (2008)

\bibitem{Mainardi:14}  Mainardi, F.: On some properties of the Mittag-Leffler function $E_\alpha(-t^\alpha)$,
completely monotone for $t>0$ with $0<\alpha<1$.  Discrete Contin.
Dyn. Syst. Ser. B. in press.

\bibitem{Metzler:00}  Metzler  R.,  Klafter, J.: The random walk's guide to anomalous
diffusion: A fractional dynamics approach.  Phys. Rep. {\bf339},
 1-77 (2000)

\bibitem{Mustapha:11}
 Mustapha, K.,  Brunner, H., Mustapha, H.,  Sch\"otzau, D.: An
hp-version discontinuous Galerkin method for integro-differential
equations of parabolic type.  SIAM J. Numer. Anal. {\bf49},
1369-1396 (2011)

\bibitem{Mustapha:13}
Mustapha, K.A.: A Superconvergent discontinuous Galerkin method for
Volterra integro-differential equations, Mathematics of Computation.
{\bf82}, 1987-2005 (2013)

\bibitem{Podlubny:99} Podlubny, I.:  Fractional Differential Equations. Academic Press, New York, USA
(1999)

\bibitem{Seybold:08}  Seybold, D.,  Hilfer, R.: Numerical algorithm for calculating the generalized Mittag-Leffler
Function. SIAM J. Numer. Anal. {\bf47},  69-88 (2008)

\bibitem{Schoetzau:00}
Sch\"otzau, D., Schwab, C.:  An hp a priori error analysis of the DG
time-stepping method for initial value problems. Calcolo, {\bf37},
207-232 (2000)


\end{thebibliography}
\end{document}